\documentclass[letter,10.5pt]{article}%
\usepackage[latin9]{inputenc}
\usepackage{geometry}
\usepackage{color}
\usepackage{amsmath}
\usepackage{amsthm}
\usepackage{amssymb}
\usepackage{graphicx}
\usepackage{setspace}
\usepackage{esint}
\usepackage{amsfonts}
\usepackage{comment}
\usepackage{booktabs}
\usepackage{lscape}
\usepackage{chngpage}
\usepackage{array}
\usepackage[titletoc,title]{appendix}
\usepackage[labelsep=period]{caption}
\usepackage{epstopdf}
\usepackage{natbib}
\usepackage[colorlinks,citecolor=blue,urlcolor=blue,filecolor=blue]{hyperref}
\usepackage{diagbox}
\usepackage{mathtools}
\usepackage{threeparttable}
\usepackage{multirow}
\usepackage{float}
\usepackage{bm}
\usepackage{siunitx}
\usepackage{natbib}%
\providecommand{\U}[1]{\protect\rule{.1in}{.1in}}

\providecommand{\U}[1]{\protect\rule{.1in}{.1in}}
\providecommand{\U}[1]{\protect\rule{.1in}{.1in}}

\numberwithin{equation}{section}
\makeatletter
\theoremstyle{plain}

\theoremstyle{definition}

\theoremstyle{plain}

\RequirePackage[colorlinks,citecolor=blue,urlcolor=blue]{hyperref}
\RequirePackage{hypernat}
\providecommand{\U}[1]{\protect\rule{.1in}{.1in}}
\newtheorem{lemma}{Lemma}
\newtheorem{theorem}{Theorem}

\newtheorem{proposition}{Proposition}
\newtheorem{definition}{Definition}
\newtheorem{example}{Example}
\newtheorem{assumption}{Assumption}
\newtheorem{remark}{Remark}

\providecommand{\definitionname}{Definition}
\providecommand{\propositionname}{Proposition}
\providecommand{\theoremname}{Theorem}
\providecommand{\definitionname}{Definition}
\providecommand{\propositionname}{Proposition}
\providecommand{\theoremname}{Theorem}
\providecommand{\definitionname}{Definition}
\providecommand{\propositionname}{Proposition}
\providecommand{\theoremname}{Theorem}
\makeatother
\providecommand{\definitionname}{Definition}
\providecommand{\propositionname}{Proposition}
\providecommand{\theoremname}{Theorem}

\graphicspath{{figures/}}
\bibpunct[, ]{(}{)}{,}{a}{}{,}

\begin{document}

\title{Extended Dynamic Programming Principle and Applications to
Time-Inconsistent Control}
\author{Yuhong Xu\thanks{Center for Financial Engineering and Math Center for
Interdiscipline Research, School of Mathematical Sciences, Soochow University,
P. R. China. The work of the first author is supported by Natural Science
Foundation of China (No.12271391; No.11871050)}, \ \ \ Shuzhen
Yang\thanks{Shandong University-Zhong Tai Securities Institute for Financial
Studies, Shandong University, PR China, (yangsz@sdu.edu.cn). This work of the
second author is was supported by National Natural Science Foundation of China (Grant No.12471450), and Taishan Scholar Talent Project Youth Project.}}
\maketitle
\date{}

\textbf{Abstract}: Since \citet{P93} established a local maximum principle for a general stochastic control problem governed by forward-backward stochastic differential equations (FBSDEs), the corresponding partial differential equation (PDE) characterization has not been developed yet. The main difficulty stems from the potential time inconsistency inherent in this class of control problems. In a dimension-augmented space, we first establish an extended dynamic programming principle (DPP). Consequently, an extended Hamilton-Jacobi-Bellman (HJB) equation is derived. The existence and uniqueness of a new type of viscosity solution is also investigated for this extended HJB equation. Compared to extant research on the stochastic maximum principle, the present paper is the first normal work on the PDE method for a control system with states evolving in both forward and backward manners. Interestingly, our extended DPP provides an equilibrium solution for general time-inconsistent control problems associated with the traditional mean-variance model, risk-sensitive control and utility optimization for narrow framing investors, among others.



\section{Introduction}

In the present paper, we study the following value function of a general cost
functional,
\begin{equation}
v(t,x)=\inf_{u(\cdot)\in\mathcal{U}[t,T]}\mathbb{E}_{t}\left[  \int_{t}%
^{T}f\left(  s,X_{s}^{t,x,u},Y_{s}^{t,x,u},Z_{s}^{t,x,u},u_{s}\right)
\mathrm{d}s+G(X_{T}^{t,x,u})+\gamma(Y_{t}^{t,x,u})\right]  ,
\label{in-value-1}%
\end{equation}
where $\mathbb{E}_{t}[\cdot]$ is the conditional expectation at time $t$, $\mathcal{U}[t,T]$ denotes all the adaptive and square integral stochastic processes on $U\in\mathbb{R}$ and
the controlled states $((X_{s}^{t,x,u}),(Y_{s}^{t,x,u},Z_{s}^{t,x,u}))$ are
solutions to the following stochastic differential equation (SDE):
\begin{equation}
X_{s}^{t,x,u}=x+\int_{t}^{s}\mu(r,X_{r}^{t,x,u},u_{r})\mathrm{d}r+\int_{t}%
^{s}\sigma(r,X_{r}^{t,x,u},u_{r})\mathrm{d}B_{r}, \label{state-1}%
\end{equation}
and backward stochastic differential equation (BSDE):
\begin{equation}
Y_{s}^{t,x,u}=\Phi\left(  X_{T}^{t,x,u}\right)  +\int_{s}^{T}h(r,X_{r}%
^{t,x,u},Y_{r}^{t,x,u},Z_{r}^{t,x,u},u_{r})\mathrm{d}r-\int_{s}^{T}%
Z_{r}^{t,x,u}\mathrm{d}B_{r}\text{, \ \ \ } \label{state2}%
\end{equation}
with $f,G,\gamma,\mu,\sigma,\Phi$ and $h$ being given deterministic
functions, and $(u_{r})_{t\leq r\leq T}$ a given control process.

\citet{P93} was the first to study the {control problem \eqref{in-value-1}}$\thicksim
${\eqref{state2}}, establishing a local maximum principle. Notably, when
$\gamma(\cdot)$ is a nonlinear function, the term $\gamma(Y_{t}^{t,x,u})$ may
cause time inconsistency in the value function $v(t,x)$. That is, the traditional DPP does not hold true. The {control problem \eqref{in-value-1}}$\thicksim$\eqref{state2} includes several interesting
examples of time-inconsistent decisions, such as those associated with the classical mean-variance
model, risk-sensitive control under a nonlinear expectation, and utility maximization for a narrow
framing investor. See Section \ref{sec:time-exam} for more details regarding
applications. To define an ``optimal" solution for a time-inconsistent control problem, the usual approach is to assume that the player and their future selves optimize simultaneously. If there exists a Nash equilibrium over time, we refer to it as the time-consistent equilibrium solution for this problem. See, for example, \citet{Basak2010} or \citet{BM2010}. So a natural question arises: does the above problem have an equilibrium solution as in \citet{BM2010}, and how can we find this time-consistent equilibrium?

When $\gamma(\cdot)$ has continuous second-order derivatives in $\mathbb{R}$,
we can apply It\^{o}'s formula to $\gamma(Y_{t}^{t,x,u})$ and transform problem \eqref{in-value-1}  into a \emph{time-consistent} problem,
\begin{equation}
v(t,x)=\inf_{u(\cdot)\in\mathcal{U}[t,T]}\mathbb{E}_{t}\left[  \int_{t}%
^{T}f\left(  s,X_{s}^{t,x,u},Y_{s}^{t,x,u},Z_{s}^{t,x,u},u_{s}\right)
\mathrm{d}s+G(X_{T}^{t,x,u})\right]  . \label{in-value-3}%
\end{equation}
With an abuse of notation, the running and terminal cost functions are still denoted
as $f$ and $G$. Thus a DPP can be established for {\eqref{in-value-3}}. However,
this DPP differs from the traditional one in that it is for an extended
control system. Specifically, we also require an auxiliary DPP for $Y^{t,x,u^{\ast}}%
$, where $u^{\ast}$ is an optimal feedback control in \eqref{in-value-3}. See Section
\ref{sec:General_framework} for more details. In the following section, we investigate the
problem \eqref{in-value-3}. Then we return to problem \eqref{in-value-1} and
its applications in Sections \ref{sec:time-in} and \ref{sec:time-exam}, respectively.
We show that the \emph{optimal control} of problem \eqref{in-value-3} provides an \emph{equilibrium }solution to problem \eqref{in-value-1}.

This paper is closely related to that of \citet{Basak2010}, the
first to employ dynamic programming while studying optimal dynamic
mean-variance problem. Using the law of total variance, they obtained a recursive
formulation for the mean-variance value function; this allowed the use of
dynamic programming in the derivation of its associated PDE. To the best of our knowledge, for the above control problem
\eqref{in-value-1}$\thicksim${\eqref{state2}, due to its time-inconsistency, prior studies in the field of stochastic control have not yet constructed an applicable DPP
and derived the related state-dependent HJB equation.} Another milestone work is \citet{BM2010} who used
a verification theorem to study time-inconsistent control problem in an extensive framework of mean-variance type,
when considering that the PDE system has smooth solutions. However, their work lacks a rigorous derivation of the associated PDE as well as an analysis of its viscosity solutions. Differing from \citet{BM2010}, from the very beginning,  we establish an extended DPP for the value function (\ref{in-value-3}), which allows a rigourous derivation of the extended HJB equation. Furthmore, we investigate the uniqueness and existence of
a viscosity solution to the extended HJB equation. Based on our DPP, we show that
the value function of certain time-inconsistent optimal control problem
developed by \citet{BM2010} is exactly a unique viscosity solution to an
extended HJB equation. Accordingly, we solve an open problem proposed by \citet{BM2010} who claimed that the theory of its viscosity solution is quite technically complicated. Since the extended HJB equation is a multidimensional PDE, for which comparison principles typically fail to apply, the analysis of its viscosity solutions has long been considered a challenging problem.

The traditional optimal control problem considers a particular case of
cost functional (\ref{in-value-3}), where $f$ does not depend on the backward states
$(Y_{s}^{t,x,u},Z_{s}^{t,x,u})_{t\leq s\leq T}$. It is well-known that two key methods, the stochastic maximum principle and the DPP, are used to
study optimal control problems. For the traditional optimal control problem,
under mild conditions, one can prove that the value function $v$ is a unique
viscosity solution of an HJB equation. We refer readers to \citet{P90,P92} and
the monographs \citet{FM06} and \citet{Y99} for the basic theory of the maximum
principle and viscosity solution of HJB equation regarding the optimal control
problem. Subsequent developments for recursive optimal control are referred to \citet{P97}, \citet{WY08}, \citet{BL08}, \citet{WY10}, \citet{LW14}, \citet{Q18}, \citet{M20}. {It is worth noting that, in classical recursive stochastic control problems, the BSDE typically represents the cost functional, while the SDE governs the state dynamics. In contrast, in our setting, both the BSDE and SDE jointly constitute the state process. When the BSDE (or expectation) becomes part of the state, the control problem transforms into a form of mean-field type, which is likely to induce time inconsistency.}

For the stochastic maximum principle of the value function (\ref{in-value-3})
when $f$ depends on the states $(Y_{s}^{t,x,u},Z_{s}^{t,x,u})$, the first
results were established by \citet{P93}. Subsequent results were developed
by \citet{Ba08}, \citet{M09}, \citet{W13}, \citet{Y10}, \citet{H17}, and
\citet{HJX19}, among others. To the best of our knowledge, no prior study has worked on either DPP or state-dependent PDE method for the value functions
(\ref{in-value-3}) or {\eqref{in-value-1}}. This difficulty arises from the fact that the
cost function depends not only  on the forward diffusion $(X^{t,x,u})$ but also
on the solution $(Y^{t,x,u},Z^{t,x,u})$ of a backward equation. In
sprit of \citet{Basak2010} and \citet{BM2010}, we predict that this control
system should be associated with a vector-valued PDE called the
extended HJB equation. When the BSDE \eqref{state2} in the state degenerates
to the classical linear expectation, our control system covers the one
investigated by \citet{BM2010}, who mainly studies the general problem of mean-variance type. See \citet{BMZ2014} and \citet{D20} for further developments.
The main contribution of this paper is to establish a DPP for the value function (\ref{in-value-3}) and to rigorously derive the associated extended HJB equation. Both the smooth and viscosity solutions are investigated
for the extended HJB equation. Hence, our paper is the first normal work that provides a strict PDE method for a
controlled system with FBSDE states.

To derive the DPP and develop the related
extended HJB equation to \eqref{state2} and \eqref{in-value-3}, we define the
auxiliary value function
\[
g(t,x)=Y_{t}^{t,x,u^*},
\]
where $u^{\ast
}(\cdot)$ is an optimal feedback  control of value function (\ref{in-value-3}).

The main contributions of this study are fourfold:

(\textbf{i}) An extended DPP is established for the non-traditional problem (\ref{in-value-3}).  Subsequently, we formally derive an extended HJB equation:%
\begin{equation}%
\begin{cases}
\underset{u\in{U}}{\inf}\{\mathcal{D}{v}\left(  t,x\right)
+f\left(  t,x,g(t,x),\sigma^{\top}\partial_{x}g(t,x),u\right)
\}=0, & \\
\text{ \ \ \ \ \ }\mathcal{D}^{\ast}g\left(  t,x\right)  +h\left(
t,x,g(t,x),\sigma^{\ast\top}\partial_{x}g(t,x),u^{\ast}(t,x)\right)  =0, &
\end{cases}
\label{ehjb}%
\end{equation}
where $v(T,x)=G(x)$,
and $g(T,x)=\Phi\left(  x\right)  $. The infinitesimal operator $\mathcal{D}$
is $\mathcal{D}(\cdot)=[\partial_{t}+\mu\partial_{x}+\frac{1}{2}\sigma
\sigma^{\top}\partial_{xx}](\cdot)$, $\mathcal{D}^{\ast}(\cdot)=[\partial
_{t}+\mu^{\ast}\partial_{x}+\frac{1}{2}\sigma^{\ast}\sigma^{\ast\top}%
\partial_{xx}](\cdot)$, $\mu^{\ast}(t,x)=\mu(t,x,u^{\ast}(t,x))$ and
$\sigma^{\ast}(t,x)=\sigma(t,x,u^{\ast}(t,x))$.

(\textbf{ii}) When the value functions $v(t,x),g(t,x)\in C^{1,2}([0,T]\times\mathbb{R}^{d})$, we
verify that $v(t,x)$ and $g(t,x)$ satisfy the extended HJB equation (\ref{ehjb}). Thus, the
verification theorem is given for this general problem.

(\textbf{iii}) Based on a modified definition of viscosity solution, we show that $v(t,x)$ and $g(t,x)$ are the unique viscosity solution of the extended HJB equation. This fills in the gap left by \citet{BM2010}.

(\textbf{iv})  We address the time-consistent solution provided by (\ref{in-value-3}) is, in fact, an equilibrium defined in \citet{BM2010}, to problem (\ref{in-value-1}). Several typical applications to time-inconsistent control problems in finance and insurance are investigated.

Note that the extended HJB equation \eqref{ehjb} is a vector valued PDE. The
theory of viscosity solution for vector-valued PDEs is not well-developed due
to the absence of a general comparison theorem. If $f\left(  t,x,y,z,u\right)  $
does not depend on $z$, then the comparison theorem holds for equation
\eqref{ehjb}. Hence, we can define the viscosity solution in the traditional
way and further prove its uniqueness. However, when $f\left(
t,x,y,z,u\right)  $ depends on $z$, the comparison theorem does not hold in
general. See \citet{xu2016} for a detailed description in the language of BSDE. To solve the problem,
we define a new version of viscosity solution by imposing additional
first-order smoothness to the coefficients of the state equations. In fact, there is little
difference between the differential functions of the first-order and
Lipschitz-continuous functions in that Lipschitz-continuity implies
differentiability almost everywhere. Our results on viscosity solution complement the existing literature on viscosity solutions to multidimensional PDEs.

Additionally,
our method of extended DPP could provide equilibrium solutions for a
general type of time-inconsistent control problem. Therefore, we provide three
examples to verify our main results, including the classical mean-variance
model, utility optimization for a narrow framing investor and risk-sensitive
portfolio optimization under a nonlinear expectation. All the three
examples can be transformed into our model and thus have an equilibrium solution that supports our method.

The remainder of this paper is organized as follows: In Section \ref{sec:General_framework}, we
formulate a general optimal control problem, in which the cost functional is
driven by a coupled FBSDE. The DPP for the optimal control problem is
investigated, and a verification theorem is established. In Section
\ref{sec5}, we introduce a new type of viscosity solution and prove that the value
function is the unique viscosity solution of the extended HJB equation. In
Section \ref{sec:time-in}, we establish the connection between time-inconsistent
control problems and our general optimal control problem. As applications, we
consider some interesting examples in Section \ref{sec:time-exam}. Finally,  Section \ref{sec8} concludes the paper.

\section{A Time-Consistent But Non-traditional Framework\label{sec:General_framework}}

Let $t$ be the starting time, and $(B_s)_{s\geq t}$ be an $m$-dimensional standard Brownian motion defined on a complete
filtered probability space $(\Omega,\mathcal{F},P;\{ \mathcal{F}_{s}\}_{s\geq
t})$, where $\{ \mathcal{F}_{s}\}_{s\geq t}$ is the $P$-augmentation of the
natural filtration generated by the Brownian motion $B$. Let $T>0$ be given.
The controlled state processes $\left(  X^{t,x,u}_{s}\right)  _{s\in\left[
t,T\right]  }$ and $\left(  Y^{t,x,u}_{s}\text{, }Z^{t,x,u}_{s}\right)
_{s\in\left[  t,T\right]  }$ are given by
\begin{equation}
\mathrm{d}X^{t,x,u}_{s}=\mu(s,X^{t,x,u}_{s},u(s))\mathrm{d}s
+\sigma(s,X^{t,x,u}_{s},u(s))\mathrm{d}B_{s}\text{,
\ \ }X^{t,x,u}_{t}=x\in\mathbb{R}\text{,} \label{stateSDE}%
\end{equation}
and
\begin{equation}
\mathrm{d}Y^{t,x,u}_{s}=-h(s,X^{t,x,u}_{s},Y^{t,x,u}_{s},Z^{t,x,u}%
_{s},u(s))\mathrm{d}s+Z^{t,x,u}_{s}\mathrm{d}B_{s}\text{,
\ \ \ } Y^{t,x,u}_{T}=\Phi\left(  X^{t,x,u}_{T}\right)  , \label{stateBSDE}%
\end{equation}
where the coefficients are functions $\mu:[0,T]\times\mathbb{R}^{d}%
\times{U}\rightarrow\mathbb{R}^{d}$, $\sigma:[0,T]\times\mathbb{R}%
^{d}\times{U}\rightarrow\mathbb{R}^{d\times m}$, $h:[0,T]\times
\mathbb{R}^{d}\times\mathbb{R}\times\mathbb{R}^{m}\times{U}%
\rightarrow\mathbb{R}$, $\Phi:\mathbb{R}^{d}\rightarrow\mathbb{R}$. We restrict ourself to the control set $\mathcal{U}[t,T]=$\{$u(\cdot
):[t,T]\times\Omega\rightarrow {U},\ \{\mathcal{F}_{s}\}\text{-adapted,} \ \mathbb{E}_t[\int_t^Tu^2(s)ds]<+\infty$\}, where ${U}$ is a subset of $\mathbb{R}^{k}$ with a given positive integer $k$.

For given deterministic functions $f$ and $G$, we introduce the following cost
functional:
\begin{equation}
J(t,x;u(\cdot))=\mathbb{E}_{t}\left[  \int_{t}^{T}f\left(  s,X^{t,x,u}%
_{s},Y^{t,x,u}_{s},Z^{t,x,u}_{s},u(s)\right)  \mathrm{d}%
s+G(X_{T}^{t,x,u})\right]  . \label{cost-0}%
\end{equation}
Then, the value function is defined as
\begin{equation}
\label{value-v}v(t,x)=\inf_{u(\cdot)\in\mathcal{U}[t,T]}J(t,x,u(\cdot)).
\end{equation}
\emph{Note that the above control problem is time-consistent} but non-traditional because its DPP is not one-dimensional; it also depends on an auxiliary dimension generated by the solution $Y^{t,x,u}_{t}$ of a BSDE. So far, no literature has provided the DPP for this system.

For simplicity of notation, we omit the time variable in $\mu,\sigma
,h,f$. We assume that $f,G$ are uniformly continuous and of polynomial growth
on independent variables, and $\mu,\sigma,h,\Phi$ satisfy the following
Lipschitz and linear growth conditions.

\begin{assumption}
\label{ass-1}There exists a constant $c>0$ such that%
\[%
\begin{array}
[c]{ll}
& \left\vert \mu(x_{1},u_{1})-\mu(x_{2},u_{2})\right\vert +\left\vert
\sigma(x_{1},u_{1})-\sigma(x_{2},u_{2})\right\vert \\
& +\left\vert h(x_{1},y_{1},z_{1},u_{1})-h(x_{2},y_{2},z_{2},u_{2})\right\vert
+\left\vert \Phi(x_{1})-\Phi(x_{2})\right\vert \\
\leq & c\left(  \left\vert x_{1}-x_{2}\right\vert +\left\vert y_{1}%
-y_{2}\right\vert +\left\vert z_{1}-z_{2}\right\vert +\left\vert u_{1}%
-u_{2}\right\vert \right)  ,
\end{array}
\]
$\forall(x_{1},y_{1},z_{1},u_{1}),(x_{2},y_{2},z_{2},u_{2})\in{\mathbb{R}^{d}%
}\times\mathbb{R}\times\mathbb{R}^{m}\times U$.
\end{assumption}

\begin{assumption}
\label{ass-2} There exists a constant $c>0$ such that
\[
\left\vert \mu(x,u)\right\vert +\left\vert \sigma(x,u)\right\vert +\left\vert
\Phi(x)\right\vert \leq c(1+\mid x\mid+\mid u\mid)\text{,}
\]
and
\[
\left\vert h(x,y,z,u)\right\vert \leq c(1+\mid x\mid+\mid y\mid+\mid
z\mid+\mid u\mid)\text{,} \quad\forall(x,y,z,u)\in{\mathbb{R}^{d}}%
\times\mathbb{R}\times\mathbb{R}^{m}\times U.
\]

\end{assumption}

Let Assumptions \ref{ass-1} and \ref{ass-2} hold. By \citet{Y99}, for
any given $(t,x)\in\lbrack0,T)\times\mathbb{R}^{d}$, FBSDE (\ref{stateSDE}%
)--(\ref{stateBSDE}) admits a unique solution $\left(  X_{s}^{t,x,u}\right)
_{s\in\left[  t,T\right]  }$ and $\left(  Y_{s}^{t,x,u}\text{, }Z_{s}%
^{t,x,u}\right)  _{s\in\left[  t,T\right]  }$ in square-integrable spaces.

{To establish the DPP for the value function $v(t,x)$, we need to
consider a weak control framework $(\Omega,\mathcal{F},P;\{\mathcal{F}%
_{s}\}_{s\geq t},B,u)$. For convenience, we still use the notation
$\mathcal{U}[t,T]$ to denote $(\Omega,\mathcal{F},P;\{\mathcal{F}_{s}\}_{s\geq
t},B,u)$. In the following, we assume that the value function (\ref{value-v}) admits an optimal feedback control $u^*(\cdot)$. In fact, this can be proved  by applying the method of $\mathcal{S}$-topology
developed in \citet{Ja96} and \citet{Ja97}. For instance, under certain conditions for the coefficients
$\mu,\sigma$ in (\ref{state-1}) and assuming $f$ does not depend on the states
$(Y^{t,x,u},Z^{t,x,u})$, Theorem 6.4 in Chapter VI of \citet{FR75} shows that
the optimal control of the value function is Lipschitz-continuous in variable
$x$. Based on the derived extended HJB equation, we can determine the optimal control and verify whether it is of feedback form. } We will first establish a DPP for $v$ and $g$. Then the associated extended HJB equation is derived. A verification theorem
is given at last.

\subsection{Dynamic programming principle}

Using the same method given in
Lemma 3.2 of \citet{Y99}, we introduce the following lemma.

\begin{lemma}
\label{DPP-2} Let $({t,x})\in[0,T]\times\mathbb{R}^{d}$. Then, for any given
$s\in\lbrack t,T] $ and $(\mathcal{F}_{s})_{s\in[t,T]}$-progressively
measurable process $\xi(\cdot)$, we have
\begin{equation}%
\begin{array}
[c]{c}%
J(s,\xi_{s};u(\cdot))=\displaystyle\mathbb{E}_{s}\bigg{[}\int_{s}%
^{T}f(r,\Theta^{s,\xi_{s},u}_{r},u(r))\mathrm{d}r+G
(X^{s,\xi_{s},u}_{T})\bigg{]},
\end{array}
\label{DPP_1-2}%
\end{equation}
where $\Theta^{s,\xi_{s},u}_{r}=(X^{s,\xi_{s},u}_{r}, Y^{s,\xi_{s},u}%
_{r},Z^{s,\xi_{s},u}_{r}),\ s\leq r\leq T$.
\end{lemma}

Now, we investigate the DPP for the value function $v(t,x)$.

\begin{theorem}
\label{DPP_1}Let Assumptions \ref{ass-1} and \ref{ass-2} hold. Then, for any
given $(t,x)\in\lbrack0,T)\times\mathbb{R}^{d}$, and $s\in\lbrack t,T]$,
\begin{equation}%
\begin{array}
[c]{l}%
v({t,x})=\displaystyle\inf_{u(\cdot)\in\mathcal{U}[t,s]}\mathbb{E}%
_{t}\bigg{[}\int_{t}^{s}f(r,\Theta_{r}^{t,x,{u}},u(r
))\mathrm{d}r+v(s,X_{s}^{{t,x},{u}})\bigg{]},\\
\end{array}
\label{DPP-1}%
\end{equation}
where $\Theta_{r}^{t,x,{u}}=(X_{r}^{t,x,u},Y_{r}^{t,x,u},Z_{r}^{t,x,u}%
),\ t\leq r\leq s$, is the solution of the following FBSDE:%
\[%
\begin{cases}
\mathrm{d}X_{r}^{t,x,u}=\mu(r,X_{r}^{t,x,u},u(r))\mathrm{d}%
r+\sigma(r,X_{r}^{t,x,u},u(r))\mathrm{d}B_{r}, & X_{t}%
^{t,x,u}=x\text{,}\\
\mathrm{d}Y_{r}^{t,x,u}=-h(r,X_{r}^{t,x,u},Y_{r}^{t,x,u},Z_{r}^{t,x,u}%
,u(r))\mathrm{d}r+Z_{r}^{t,x,u}\mathrm{d}B_{r}, & Y_{s}%
^{t,x,u}=Y_{s}^{t,x,u^{\ast}},
\end{cases}
\]
and $\Theta_{r}^{t,x,{u}^{\ast}}=(X_{r}^{t,x,u^{\ast}},Y_{r}^{t,x,u^{\ast}%
},Z_{r}^{t,x,u^{\ast}}),\ s\leq r\leq T$ is the solution of the following
FBSDE:%
\[%
\begin{cases}
\mathrm{d}X_{r}^{t,x,u^{\ast}}=\mu(r,X_{r}^{t,x,u^{\ast}},u^{\ast}%
(r,X_{r}^{t,x,u^{\ast}}))\mathrm{d}r+\sigma(r,X_{r}^{t,x,u^{\ast}},u^{\ast
}(r,X_{r}^{t,x,u^{\ast}}))\mathrm{d}B_{r}\text{, }X_{s}^{t,x,u^{\ast}}%
=X_{s}^{t,x,u}, & \\
\mathrm{d}Y_{r}^{t,x,u^{\ast}}=-h(r,X_{r}^{t,x,u^{\ast}},Y_{r}^{t,x,u^{\ast}%
},Z_{r}^{t,x,u^{\ast}},u^{\ast}(r,X_{r}^{t,x,u^{\ast}}))\mathrm{d}%
r+Z_{r}^{t,x,u^{\ast}}\mathrm{d}B_{r}\text{, }Y_{T}^{t,x,u^{\ast}}=\Phi\left(
X_{T}^{t,x,u^{\ast}}\right).&
\end{cases}
\]
where $u^{\ast}(\cdot)$ denotes an optimal feedback control of the value function
(\ref{value-v}).
\end{theorem}

\noindent\textbf{Proof.} We denote the right side of equation (\ref{DPP-1}) by
$\tilde{v}(t,x)$. For any given sufficiently small $\varepsilon>0$, based on
the definition of $v({t,x})$, there exists
\[
\displaystyle{u}^{s-t}(r)=%
\begin{cases}
u_{1}(r), & \mbox{ $t\leq r< s$},\\
{u}^{\ast}(r,x), & \mbox{ $ s\leq r\leq T$},
\end{cases}
\]
which belongs to $\mathcal{U}[t,T]$ such that%

\[%
\begin{array}
[c]{cl}
& v({t,x})+\varepsilon\\
> & J(t,x;u^{s-t}(\cdot))\\
= & \displaystyle\mathbb{E}_{t}\bigg{[}\int_{t}^{T}f(r,\Theta_{r}%
^{{t,x},u^{s-t}},u^{s-t}(r))\mathrm{d}r+G(X_{T}%
^{t,x,u^{s-t}})\bigg{]}\\
= & \displaystyle\mathbb{E}_{t}\bigg{[}\int_{t}^{s}f(r,\Theta_{r}%
^{{t,x},u^{s-t}},u^{s-t}(r))\mathrm{d}r\\
& \displaystyle+\mathbb{E}_{t}[\int_{s}^{T}f(r,\Theta_{r}^{{t,x},u^{s-t}%
},u^{s-t}(r))\mathrm{d}r+G(X_{T}^{t,x,u^{s-t}}%
)\mid\mathcal{F}_{s}]\bigg{]}\\
= & \displaystyle\mathbb{E}_{t}\bigg{[}\int_{t}^{s}f(r,\Theta_{r}^{{t,x}%
,u_{1}},u_{1}(r))\mathrm{d}r\\
& \displaystyle+\mathbb{E}_{s}[\int_{s}^{T}f(r,\Theta_{r}^{s,X_{s}^{t,x,u_{1}%
},{u}^{\ast}},{u}^{\ast}(r,X_{r}^{{t,x},u_{1}}))\mathrm{d}r+G(X_{T}%
^{s,X_{s}^{t,x,u_{1}},{u}^{\ast}})\bigg{]}\\
= & \displaystyle\mathbb{E}_{t}\bigg{[}\int_{t}^{s}f(r,\Theta_{r}^{{t,x}%
,u_{1}},u_{1}(r))\mathrm{d}r+J(s,X_{s}^{{t,x},u_{1}}%
;{u}^{\ast}(\cdot))\bigg{]}\\
\geq & \displaystyle\mathbb{E}_{t}\bigg{[}\int_{t}^{s}f(r,\Theta_{r}%
^{{t,x},u_{1}},u_{1}(r))\mathrm{d}r+v(s,X_{s}^{{t,x}%
,u_{1}})\bigg{]}\\
\geq & \tilde{v}({t,x}).
\end{array}
\]
The second equality is derived from Lemma \ref{DPP-2}. Conversely, by the
definition of the value function $v(t,x)$, for a given
\[
\displaystyle u^{s-t}(r)=%
\begin{cases}
u_{1}(r), & \mbox{ $t\leq r< s$},\\
{u}^{\ast}(r,x), & \mbox{ $ s\leq  r\leq T$},
\end{cases}
\]
we have
\[
v({s,X_{s}^{t,x,u_{1}}})=J(s,X_{s}^{t,x,u_{1}};u^{\ast}(\cdot)).
\]
Note that $v({t,x})\leq J(t,x;u^{s-t}(\cdot))$,
which implies that
\[
v({t,x})\leq\displaystyle\mathbb{E}_{t}\bigg{[}\int_{t}^{s}f(r,\Theta
_{r}^{{t,x},u_{1}},u_{1}(r))\mathrm{d}r+J(s,X_{s}%
^{{t,x},u_{1}};{u}^{\ast}(\cdot))\bigg{]}.
\]
Therefore, $v({t,x})\leq\tilde{v}({t,x}).$
This completes the proof. $\ \ \ \ \ \ \ \Box$

\subsection{Extended HJB equation}

To investigate the related HJB equation for the value function $v(t,x)$ in
(\ref{value-v}). We introduce another value function for the solution of BSDE
(\ref{stateBSDE}),
\begin{equation}
g(t,x)=Y_{t}^{t,x,u^*}, \label{value-g}%
\end{equation}
where $u^{\ast}(\cdot)$ denotes an optimal feedback control of the value function
(\ref{value-v}).
Then the related
extended HJB equation is given as follows:%
\begin{align}
\underset{u\in U}{\inf}\{\mathcal{D}{v}\left(  t,x\right)
+f\left(  t,x,g(t,x),\sigma^{\top}\partial_{x}g(t,x),u\right)  \}  &
=0,\label{hjb-1}\\
\mathcal{D}^{\ast}g\left(  t,x\right)  +h\left(  t,x,g(t,x),\sigma^{\ast\top
}\partial_{x}g(t,x),u^{\ast}(t,x)\right)   &  =0, \label{hjb-10}%
\end{align}
subject to $v(T,x)=G(x)$ and $g(T,x)=\Phi\left(  x\right)$. The infinitesimal
operator $\mathcal{D}$ is $\mathcal{D}(\cdot)=[\partial_{t}+\mu\partial
_{x}+\frac{1}{2}\sigma\sigma^{\top}\partial_{xx}](\cdot)$, $\mathcal{D}^{\ast
}(\cdot)=[\partial_{t}+\mu^{\ast}\partial_{x}+\frac{1}{2}\sigma^{\ast}%
\sigma^{\ast\top}\partial_{xx}](\cdot)$ and $\mu^{\ast}(t,x)=\mu(t,x,u^{\ast
}(t,x)),\ \sigma^{\ast}(t,x)=\sigma(t,x,u^{\ast}(t,x))$.

Let $C^{1,2}([0,T]\times\mathbb{R}^{d})$ be the space of functions with
continuous first-order derivatives on $t\in[0,T]$ and continuous second-order
derivatives on $x\in\mathbb{R}^{d}$. We first establish the following result.

\begin{theorem}
\label{the-class} Suppose Assumptions \ref{ass-1}, \ref{ass-2} hold, and the value
functions $v,g\in C^{1,2}([0,T]\times\mathbb{R}^{d})$. Then $(v,g)$ is a
solution of the extended HJB equation (\ref{hjb-1})--(\ref{hjb-10}).
\end{theorem}

\noindent\textbf{Proof.} For any given $({t,x})\in\lbrack0,T]\times
\mathbb{R}^{d}$ and a control ${u}(\cdot)\in\mathcal{U}[t,s]$, define
\[
\displaystyle u_{1}(r)=%
\begin{cases}
{u}(r), & \mbox{ $t\leq r< s$},\\
{u}^{\ast}(r,x), & \mbox{ $ s\leq r\leq T$},
\end{cases}
\]
where ${u}^{\ast}(r,x)$ is an optimal control of value function $v$ on
$(s,T]$. Let $X^{t,x,u_{1}}(\cdot)$ be the corresponding solution of equation
(\ref{stateSDE}), and $(Y^{t,x,u_{1}}(\cdot),Z^{t,x,u_{1}}(\cdot))$ be the
solution of equation (\ref{stateBSDE}). According to Theorem \ref{DPP_1}, we
have
\[
v({t,x})\leq\displaystyle\mathbb{E}_{t}\bigg{[}\int_{t}^{s}f(r,\Theta
_{r}^{{t,x},u},u(r))\mathrm{d}r+v(s,X_{s}^{{t,x},{u}})\bigg{]},
\]
which implies that
\[%
\begin{array}
[c]{cl}%
0 & \leq\displaystyle\frac{\mathbb{E}_{t}\big{[}v(s,X_{s}^{{t,x},u}%
)-v({t,x})\big{]}}{s-t}+\frac{1}{s-t}\mathbb{E}_{t}\int_{t}^{s}f(r,\Theta
_{r}^{{t,x},u},u(r))\mathrm{d}r\\
& =\displaystyle\frac{1}{s-t}\mathbb{E}_{t}\int_{t}^{s}\bigg{[}\mathcal{D}%
v(r,X_{r}^{t,x,{u}})+f(r,\Theta_{r}^{{t,x},u},u(r
))\bigg{]}\mathrm{d}r.
\end{array}
\]
{Letting $s\downarrow t$, since $v,g\in C^{1,2}([0,T]\times\mathbb{R}^{d})$, based on the Feynman-Kac formula of BSDE, we have}
\[
0\leq\mathcal{D}{v}\left(  t,x\right)  +f(t,x,g(t,x),\sigma^{\top}%
\partial_{x}g(t,x),u(t)).
\]
Thus, it follows that
\begin{equation}%
\begin{array}
[c]{c}%
0\leq\underset{u\in{U}}{\inf}\{\mathcal{D}{v}\left(
t,x\right)  +f\left(  t,x,g(t,x),\sigma^{\top}\partial_{x}%
g(t,x),u\right)  \}.
\end{array}
\label{estimate-3}%
\end{equation}

For any given $\varepsilon>0$, there exists a control
\[
\displaystyle u(r)=%
\begin{cases}
u^{\varepsilon}(r), & \mbox{ $t\leq r< s$},\\
{u}^{\ast}(r,x), & \mbox{ $ s\leq r\leq T$},
\end{cases}
\]
such that
\[
v({t,x})+\varepsilon(s-t)\geq\displaystyle\mathbb{E}_{t}\bigg{[}\int_{t}%
^{s}f(r,\Theta_{r}^{{t,x},u^{\varepsilon}},u^{\varepsilon}(r))\mathrm{d}r +v(s,X_{s}^{{t,x},u^{\varepsilon}%
})\bigg{]}.
\]
Thus, it follows from It\^{o}'s formula that,
\[%
\begin{array}
[c]{cl}%
\varepsilon & \geq\displaystyle\frac{\mathbb{E}_{t}\big{[}v(s,X_{s}%
^{{t,x},u^{\varepsilon}})-v({t,x})\big{]}}{s-t}+\frac{1}{s-t}\mathbb{E}%
_{t}\int_{t}^{s}f(r,\Theta_{r}^{{t,x},u},u^{\varepsilon}(r))\mathrm{d}r\\
& =\displaystyle\frac{1}{s-t}\mathbb{E}_{t}\int_{t}^{s}\bigg{[}\mathcal{D}%
v(r,X_{r}^{t,x,u^{\varepsilon}})+f(r,\Theta_{r}^{{t,x},u^{\varepsilon}%
},u^{\varepsilon}(r))\bigg{]}\mathrm{d}r.
\end{array}
\]
Therefore,
\begin{equation}%
\begin{array}
[c]{c}%
0\geq\underset{u\in{U}}{\inf}\{\mathcal{D}{v}\left(
t,x\right)  +f\left(  t,x,g(t,x),\sigma^{\top}\partial_{x}%
g(t,x),u\right)  \}.
\end{array}
\label{estimate-2}%
\end{equation}
Then, combining equations (\ref{estimate-3}) and (\ref{estimate-2}), the value
function $v$ is the solution of equation (\ref{hjb-1}).

Based on the nonlinear Feynman-Kac formula (see \citet{EP92}), the function
$g(t,x)=Y_{t}^{t,x,u^{\ast}}\in C^{1,2}([0,T]\times\mathbb{R}^{d})$ is the
classical solution of
\[
\mathcal{D}^{\ast}g\left(  t,x\right)  +h\left(  t,x,g(t,x),\sigma^{\ast\top
}\partial_{x}g(t,x),u^{\ast}(t,x)\right)  =0,\ g(T,x)=\Phi(x).
\]
The proof is complete. \ \ \ \ \ \ \ $\Box$

\subsection{Verification theorem}

Now, we establish the well-known verification theorem for value function
(\ref{value-v}).

\begin{theorem}
\label{the-ver} Let Assumptions \ref{ass-1}, \ref{ass-2} hold, and $(v,g)$ be
the classical solution of equations (\ref{hjb-1}) and \eqref{hjb-10}. Then, we have

(i) $v(t,x)\leq J(t,x;u(\cdot))$ for any given $u(\cdot)\in\mathcal{U}[t,T]$
and $(t,x)\in[0,T]\times\mathbb{R}^{d}$;

(ii) For any $(t,x)\in[0,T]\times\mathbb{R}^{d}$, if $u^{*}(\cdot)\in \mathcal{U}^*[t,T]$ satisfies $v(t,x)=J(t,x;u^{*}(\cdot))$,
then $u^{*}(\cdot)$ is an optimal control of value function
(\ref{value-v}).
\end{theorem}

\noindent\textbf{Proof.} (i) Note that $v,g$ are the classical solution of
equations (\ref{hjb-1}) and (\ref{hjb-10}). Thus for an optimal feedback control
$u^{*}(\cdot)$, and $(t,x)\in[0,T)\times\mathbb{R}^{d}$,
\[
0=\mathcal{D}^{*}{v}\left(  t,x\right)  +f\left(  t,x,g(t,x),\sigma^{*\top
}\partial_{x} g(t,x),u^{*}(t,x)\right)  ,
\]
{which deduces that}
\[
0= \mathcal{D}^{*}{v}\left(  s,X_{s}^{t,x,u^{*}}\right)  +f\left(
s,X_{s}^{t,x,u^{*}},g(s,X_{s}^{t,x,u^{*}}),\sigma^{*\top}\partial_{x}
g(s,X_{s}^{t,x,u^{*}}),u^{*}(s,X_{s}^{t,x,u^{*}})\right)  .
\]
Then, applying It\^{o}'s formula to ${v}( s,X_{s}^{t,x,u^{*}})$, we have
\[
v(t,x)=\displaystyle \mathbb{E}_{t}\left[  -\int_{t}^{T} \mathcal{D}^{*}{v}(
r,X_{r}^{t,x,u^{*}})\mathrm{d}r+{G}( X_{T}^{t,x,u^{*}})\right]  .
\]
Therefore
\begin{equation}
\label{ver-1}v(t,x)= \mathbb{E}_{t}\left[  \int_{t}^{T}f\left(  r,X_{r}%
^{t,x,u^{*}},g(r,X_{r}^{t,x,u^{*}}),\sigma^{*\top}\partial_{x} g(r,X_{r}%
^{t,x,u^{*}}),u^{*}(r,X_{r}^{t,x,u^{*}})\right)  \mathrm{d}r+{G}%
(X_{T}^{t,x,u^{*}})\right]  .
\end{equation}

Furthermore, we can verify that
\[
(X_{r}^{t,x,u^{\ast}},g(r,X_{r}^{t,x,u^{\ast}}),\sigma^{\ast\top}\partial
_{x}g(r,X_{r}^{t,x,u^{\ast}}))=(X_{r}^{t,x,u^{\ast}},Y_{r}^{t,x,u^{\ast}%
},Z_{r}^{t,x,u^{\ast}}),\ t\leq r\leq T
\]
is the unique solution of the following FBSDE:%
\[%
\begin{cases}
\mathrm{d}X_{r}^{t,x,u^{\ast}}=\mu(r,X_{r}^{t,x,u\ast},u^{\ast}(r,X_{r}%
^{t,x,u^{\ast}}))\mathrm{d}r+\sigma(r,X_{r}^{t,x,u^{\ast}},u^{\ast}%
(r,X_{r}^{t,x,u^{\ast}}))\mathrm{d}B_{r},\ X_{t}^{t,x,u^{\ast}}=x\text{,} & \\
\text{\ }\mathrm{d}Y_{r}^{t,x,u^{\ast}}=-h(r,X_{r}^{t,x,u},Y_{r}^{t,x,u^{\ast
}},Z_{r}^{t,x,u^{\ast}},u^{\ast}(r,X_{r}^{t,x,u^{\ast}}))\mathrm{d}%
r+Z_{r}^{t,x,u^{\ast}}\mathrm{d}B_{r},\ Y_{T}^{t,x,u}=\Phi(X_{T}^{t,x,u}). &
\end{cases}
\]
Thus, equation (\ref{ver-1}) can be rewritten as
\[
v(t,x)=\mathbb{E}_{t}\left[  \int_{t}^{T}f\left(  r,X_{r}^{t,x,u^{\ast}}%
,Y_{r}^{t,x,u^{\ast}},Z_{r}^{t,x,u^{\ast}},u^{\ast}(r,X_{r}^{t,x,u^{\ast}%
})\right)  \mathrm{d}r+{G}(X_{T}^{t,x,u^{\ast}})\right]  .
\]
From the optimality of $u^{\ast}(\cdot)$  for the value function
(\ref{value-v}), we get that,
\[
v(t,x)\leq J(t,x;u(\cdot)),
\]
for any given control $u(\cdot)\in\mathcal{U}[t,T]$.

(ii) From the proof in part (i), we have that for any $(t,x)\in
[0,T]\times\mathbb{R}^{d}$ and any $u(\cdot)\in\mathcal{U}[t,T]$, $v(t,x)\leq
J(t,x;u(\cdot))$. Thus, if $u^{*}(\cdot)\in\mathcal{U}[t,T]$ such that
$v(t,x)= J(t,x;u^{*}(\cdot))$, then $u^{*}(\cdot)$ is an optimal control for
value function (\ref{value-v}).\ \ \ \ \ \ \ $\Box$

\section{Viscosity Solution}\label{sec5}

In this section, we discuss the viscosity solution of the extended HJB
equation (\ref{hjb-1})--(\ref{hjb-10}) in two situations: $f$ is independent
of $Z^{t,x,u}$, and $f$ depends on $Z^{t,x,u}$. In the first case, based on
the method given in \citet{EPR97} and following the classical definition of
viscosity solution, we show that the value function $(v,g)$ is the unique
viscosity solution of (\ref{hjb-1}) and (\ref{hjb-10}). While $f$ depends on
$Z^{t,x,u}$, to the best of our knowledge, there is no result about the
viscosity solution of HJB equations (\ref{hjb-1}) and (\ref{hjb-10}). One of
the reasons is that the comparison theorem for the vector-valued
function $(v,g)$ does not hold in general. To overcome the difficulty, we impose
a first-order smoothness condition on the coefficients of the controlled FBSDE (\ref{stateSDE}) and (\ref{stateBSDE}).
See subsection \ref{f-z} for more details. As we know that, compared with the
Lipschitz-continuous condition, the first-order smoothness condition is
trivial since the uniformly Lipschitz continuity implies almost everywhere differentiability.

\subsection{When $f$ does not depend on $Z^{t,x,u}$}

Inspirited by \citet{EPR97}, we first give a definition of a viscosity
solution of the vector valued PDE (\ref{hjb-1}) and (\ref{hjb-10}) according
to the classical type:

\begin{definition}
\label{def-vis} Let $(w_{1},w_{2})\in C([0,T]\times\mathbb{R}^{d})$.
$(w_{1},w_{2})$ is a viscosity sub-solution of (\ref{hjb-1}) and
(\ref{hjb-10}), if $\forall(\Gamma_{1},\Gamma_{2})\in C^{1,2}([0,T]\times
\mathbb{R}^{d})$, for any $(t,x)\in\lbrack0,T]\times\mathbb{R}^{d}$ such that
$\Gamma_{i}\geq w_{i}$ and $\Gamma_{i}(t,x)=w_{i}({t,x}),\ i=1,2$, we have
\[%
\begin{array}
[c]{l}%
\displaystyle\underset{u\in{U}}{\inf}\{\mathcal{D}{\Gamma
}_{1}\left(  t,x\right)  +f\left(  t,x,\Gamma_{2}(t,x),u\right)  \}\geq0,
\end{array}
\]
and
\[%
\begin{array}
[c]{l}%
\displaystyle\mathcal{D}^{\ast}\Gamma_{2}\left(  t,x\right)  +h\left(
t,x,\Gamma_{2}(t,x),\sigma^{\ast\top}\partial_{x}\Gamma_{2}(t,x),u^{\ast
}(t,x)\right)  \geq0.
\end{array}
\]
$(w_{1},w_{2})$ is a viscosity super-solution of (\ref{hjb-1}) and
(\ref{hjb-10}), if $\forall(\Gamma_{1},\Gamma_{2})\in C^{1,2}([0,T]\times
\mathbb{R}^{d})$, for any $(t,x)\in\lbrack0,T]\times\mathbb{R}^{d}$ such that
$\Gamma_{i}\leq w_{i}$ and $\Gamma_{i}(t,x)=w_{i}({t,x}),\ i=1,2$, we have
\[%
\begin{array}
[c]{l}%
\displaystyle\underset{u\in{U}}{\inf}\{\mathcal{D}{\Gamma
}_{1}\left(  t,x\right)  +f\left(  t,x,\Gamma_{2}(t,x),u\right)  \}\leq0,
\end{array}
\]
and
\[%
\begin{array}
[c]{l}%
\displaystyle\mathcal{D}^{\ast}\Gamma_{2}\left(  t,x\right)  +h\left(
t,x,\Gamma_{2}(t,x),\sigma^{\ast\top}\partial_{x}\Gamma_{2}(t,x),u^{\ast
}(t,x)\right)  \leq0.
\end{array}
\]
$(w_{1},w_{2})$ is a viscosity solution of (\ref{hjb-1}) and (\ref{hjb-10}),
if it is both a viscosity super-solution and a viscosity sub-solution of
(\ref{hjb-1}) and (\ref{hjb-10}) over $[0,T]\times\mathbb{R}^{d}$.
\end{definition}

Before checking the viscosity solution, we first show that the value function
$v$ has the following continuous properties. Here, the function $f$ can depend
on $Z^{t,x,u}$.

\begin{lemma}
\label{lem-4} Let Assumptions \ref{ass-1} and \ref{ass-2} hold. The value
functions $v$ and $g$ defined in (\ref{value-v}) and (\ref{value-g}) are
continuous and of linear growth. Namely, there exists a constant $c>0$ such
that, for every $({t,x})\in\lbrack0,T]\times\mathbb{R}^{d}$, we have
\[%
\begin{array}
[c]{c}%
\left\vert v({t,x})\right\vert +\mid g({t,x})\mid\leq c(1+\mid x\mid),
\end{array}
\]
and $\forall\varepsilon>0$, there exists a constant $\delta>0,$ such that
$\sqrt{\mid t-s\mid+\mid x-y\mid}<\delta$, we have
\[
\mid v(t,x)-v(s,y)\mid+\mid g(t,x)-g(s,y)\mid<\varepsilon.
\]

\end{lemma}

\noindent\textbf{Proof.} Based on Assumption \ref{ass-2}, we obtain that $v$
is of linear growth. Next, we prove the continuity of $v$ at $(t,x)$: For a
given $u(\cdot)\in\mathcal{U}[t,T]$ and the initial data $(t,x)\in
[0,T]\times\mathbb{R}^{d}$ and $(s,y)\in[0,T]\times\mathbb{R}^{d}$ with $t\leq
s\leq T$, we have
\[%
\begin{array}
[c]{cl}
& J(t,x;u(\cdot))-J(s,y;u(\cdot))\\
= & \displaystyle\mathbb{E}_{t}\bigg{[}\int_{t}^{T}f(r,\Theta^{t,x,u}%
_{r},u(r))\mathrm{d}r -\int_{s}^{T}f(r,\Theta^{s,y,u}%
_{r},u(r))\mathrm{d}r \displaystyle +G (X^{t,x,u}_{T}%
)-G(X^{s,y,u}_{T})\bigg{]}\\
= & \displaystyle\mathbb{E}_{t}\bigg{[}\int_{t}^{s}f(r,\Theta^{t,x,u}%
_{r},u(r))\mathrm{d}r + \int_{s}^{T} \big{[}f(r,\Theta
^{t,x,u}_{r},u(r))-f(r,\Theta^{s,y,u}_{r},u(r))\big{]}\mathrm{d}r \bigg{]}\\
& +\displaystyle\mathbb{E}_{t}\bigg{[}G (X^{t,x,u}_{T})-G(X^{s,y,u}%
_{T})\bigg{]}.
\end{array}
\]
Denote
\[%
\begin{array}
[c]{cl}%
I_{1}= & \displaystyle\mathbb{E}_{t}\bigg{[}\int_{t}^{s}f(r,\Theta^{t,x,u}%
_{r},u(r))\mathrm{d}r\\
& \displaystyle+ \int_{s}^{T} \big{[}f(r,\Theta^{t,x,u}_{r},u(r))-f(r,\Theta^{s,y,u}_{r},u(r))\big{]}\mathrm{d}r
\bigg{]},\\
I_{2}= & \displaystyle\mathbb{E}_{t}\bigg{[}G (X^{t,x,u}_{T})-G(X^{s,y,u}%
_{T})\bigg{]}.
\end{array}
\]
For any given $\varepsilon>0$, by Assumption \ref{ass-1}, there exists
$\delta>0$ and $\sqrt{\mid t-s\mid+\mid x-y \mid} <\delta$ such that
$\displaystyle\left|  I_{1}\right|  +\left|  I_{2}\right|  <\frac{\varepsilon
}{2}$. Therefore, we have
\[
\left|  J(t,x;u(\cdot))-J(s,y;u(\cdot))\right|  \leq\frac{\varepsilon}{2},
\]
which further implies that $\left|  v(t,x)-v(s,y)\right|  \leq\frac
{\varepsilon}{2}$. Similarly, we can prove that
\[
\left|  g(t,x)-g(s,y)\right|  \leq\frac{\varepsilon}{2}.
\]

$\Box$

{ \bigskip}

%

\begin{theorem}
{ \label{the-vis-un0} Let Assumptions \ref{ass-1}, \ref{ass-2} hold. The value
function $(v,g)$ is the unique viscosity solution of the extended HJB equation
(\ref{hjb-1})--(\ref{hjb-10}). Namely, }

{ (i) For a given optimal feedback control $u^{*}(\cdot)$ of value function
(\ref{value-v}), the value function $g(t,x)$ of (\ref{value-g}) is the unique
viscosity solution of (\ref{hjb-10}); }

{ (ii) The value function $v(t,x)$ of (\ref{value-v}) is the unique viscosity
solution of (\ref{hjb-1}). }
\end{theorem}

{ \noindent\textbf{Proof:} For a given optimal feedback control $u^{*}(\cdot)$ of value
function (\ref{value-v}), following the results of \citet{P97}, we deduce that
$g(t,x)$ is a viscosity solution of equation (\ref{hjb-10}). Then, based on
the results of \citet{EPR97}, we deduce that the value function $v(t,x)$ is a viscosity solution
of (\ref{hjb-1}).

For a given optimal feedback control $u^{\ast}(\cdot)$ of
value function (\ref{value-v}), the uniqueness of viscosity solution of
(\ref{hjb-10}), subject to $g(T,x)=\Phi\left(  x\right)  $, is obtained
directly by the Theorem 8.2 of \citet{CIL92}. }

For a given $g(t,x)\in C([0,T]\times
\mathbb{R}^{d})$,  the
uniqueness of the viscosity solution of (\ref{hjb-1}) can be also obtained by the Theorem 8.2 of \citet{CIL92}.

For another $u^{\ast}(\cdot)$ and the related $g(t,x)$, repeat the above
procedure, we can still prove that the value function $v(t,x)$ of
(\ref{value-v}) is the unique viscosity solution of (\ref{hjb-1}). Hence the
uniqueness of viscosity solution to (\ref{hjb-1}) does not depend on either
the choice of an optimal control or the related $g(t,x)$. $\ \ \ \ \ \ \ \Box$

\begin{remark}
{ \label{re-un} Theorem \ref{the-vis-un0} tells us that the value function
(\ref{value-v}) is the unique viscosity solution of (\ref{hjb-1}), while there
may be several functions $g(t,x)$ for (\ref{hjb-10}) depending given optimal
controls $u^{\ast}(\cdot)$. In fact, equation (\ref{hjb-10}) is an auxiliary
equation. What we are concerned about is the the first equation (\ref{hjb-1}).
}
\end{remark}

\subsection{{When $f$ depends on $Z^{t,x,u}$}}

\label{f-z} We first give an example to show that Definition \ref{def-vis} is
not suitable for equations (\ref{hjb-1}) and (\ref{hjb-10}) when $f$ depends
on $Z^{t,x,u}$. Thus we need to redefine the viscosity solution.

\begin{example}
{ \label{ex-1} We consider a simple vector-valued ordinal differential
equations which are similar with equations (\ref{hjb-1}) and (\ref{hjb-10}):
\begin{equation}%
\begin{cases}
\displaystyle2-\left\vert \frac{\mathrm{d}u(x)}{\mathrm{d}x}\right\vert
-\frac{\mathrm{d}v(x)}{\mathrm{d}x}=0,\\
\displaystyle1-\left\vert \frac{\mathrm{d}v(x)}{\mathrm{d}x}\right\vert
=0,\ 0\leq x\leq2,\\
u(0)=u(2)=v(0)=v(2)=0.
\end{cases}
\label{ex-eq-1}%
\end{equation}
Applying Definition \ref{def-vis}, we can verify that
\[
\displaystyle\tilde{v}(x)=%
\begin{cases}
x, & \mbox{ $0\leq x\leq 1$,}\\
2-x, & \mbox{ $ 1< x\leq 2$},
\end{cases}
\]
is the unique viscosity solution of the second equation in (\ref{ex-eq-1}):
\[
\displaystyle1-\mid\frac{\mathrm{d}v(x)}{\mathrm{d}x}\mid=0,\ v(0)=v(2)=0.
\]
} { Now, we consider the first equation in (\ref{ex-eq-1}),
\begin{equation}
\label{ex-eq-2}\displaystyle 2-\mid\frac{\mathrm{d}u(x)}{\mathrm{d}x}\mid-
\frac{\mathrm{d}v(x)}{\mathrm{d}x}=0,\ 0\leq x\leq2, \ u(0)=u(2)=0.
\end{equation}
When $x\neq1$, we can show that
\[
\displaystyle \tilde{u}(x)=
\begin{cases}
x+2, & \mbox{ $0\leq x< 1$,}\\
6-3x, & \mbox{ $ 1< x\leq 2$},
\end{cases}
\]
is the viscosity solution of equation (\ref{ex-eq-2}). Furthermore, let
$\tilde{u}(1)=3$. Thus $\tilde{u}(\cdot)$ is continuous on interval $[0,2]$. }

{ Applying Definition \ref{def-vis}, we consider the point $x=1$, and
functions $\phi_{1}(\cdot),\phi_{2}(\cdot)\in C^{1}(\mathbb{R})$, which
satisfy
\[
\phi_{1}(x)\geq\tilde{u}(x),\ \phi_{1}(1)=\tilde{u}(1),\ \phi_{2}(x)\geq
\tilde{v}(x),\ \phi_{2}(1)=\tilde{v}(1).
\]
It follows that
\[
\displaystyle \frac{\mathrm{d}\phi_{1}(1)}{\mathrm{d}x}\in[-3,1],\ \frac
{\mathrm{d}\phi_{2}(1)}{\mathrm{d}x}\in[-1,1].
\]
Thus, when $\displaystyle \frac{\mathrm{d}\phi_{1}(1)}{\mathrm{d}x}%
=-3,\ \frac{\mathrm{d}\phi_{2}(1)}{\mathrm{d}x}=1$,
\[
\displaystyle 2-\mid\frac{\mathrm{d}\phi_{1}(1)}{\mathrm{d}x}\mid-
\frac{\mathrm{d}\phi_{2}(1)}{\mathrm{d}x}=-2<0;
\]
When $\displaystyle \frac{\mathrm{d}\phi_{1}(1)}{\mathrm{d}x}=0,\ \frac
{\mathrm{d}\phi_{2}(1)}{\mathrm{d}x}=0$,
\[
\displaystyle 2-\mid\frac{\mathrm{d}\phi_{1}(1)}{\mathrm{d}x}\mid-
\frac{\mathrm{d}\phi_{2}(1)}{\mathrm{d}x}=2>0.
\]
Hence, $\tilde{u}(\cdot)$ is not a viscosity solution of equation
(\ref{ex-eq-2}). }

{ In fact, we cannot find a viscosity solution for equation (\ref{ex-eq-2}).
Thus, Definition \ref{def-vis} is inapposite to equations (\ref{hjb-1}) and
(\ref{hjb-10}). The problem is that the term $\sigma^{\top}\partial_{x}
\Gamma_{2}(t,x)$ in Definition \ref{def-vis} may cause the ill-posedness of
the viscosity solution of equation (\ref{hjb-1}). }
\end{example}

{ To overcome the problem in the definition of viscosity solution in
Definition \ref{def-vis}, we introduce the following first-order smoothness
condition for the coefficients of FBSDE (\ref{stateSDE}) and (\ref{stateBSDE}%
), which is used to obtain the continuity of solution $Z_{t}^{t,x,u^{*}%
}=\sigma^{*\top}(t,x)\partial_{x}g(t,x)$ at $(t,x)\in[0,T]\times\mathbb{R}%
^{n}$. The following assumption and lemma come from \citet{MZ02}. }

\begin{assumption}
{ \label{ass-4} For any given optimal feedback control $u^{*}(\cdot)\in \mathcal{U}[t,T]$ of value function
(\ref{value-v}), $\mu^{*}\in C_{b}^{0,1}([0,T]\times\mathbb{R}^{d})$,
$\sigma^{*}\in C_{b}^{0,1}([0,T]\times\mathbb{R}^{d})$, $h^{*}\in C_{b}%
^{0,1}([0,T]\times\mathbb{R}^{d})$, $\Phi(x)\in C_{b}^{1}(\mathbb{R})$, i.e.,
the partial derivatives of $\mu^{*},\sigma^{*},h^{*},\Phi$ to $x\in
\mathbb{R}^{d}$ are continuous and uniformly bounded, where
\begin{equation}%
\begin{array}
[c]{cl}
& \mu^{*}(t,x)=\mu(t,x,u^{*}(t,x)),\ \sigma^{*}(t,x)=\sigma(t,x,u^{*}(t,x)),\\
& h^{*}(t,x,y,z)=h(t,x,y,z,u^{*}(t,x)),\ (t,x,y,z)\in[0,T]\times\mathbb{R}%
^{n}\times\mathbb{R}\times\mathbb{R}^{m}.
\end{array}
\end{equation}
}
\end{assumption}

{ Note that uniformly Lipschitz continuity implies absolute continuity, and
absolute continuity implies almost everywhere differentiability. Thus,
comparing with Lipschitz continuity, Assumption \ref{ass-4} is not a
strong condition. }

\begin{lemma}
{ \label{lem-5} Suppose Assumption \ref{ass-4} holds. Then, we have that }

{ (i) $Z_{s}^{t,x,u^{*}})=\sigma^{*\top}\partial_{x}Y_{s}^{t,x,u^{*}}$, where
$(Y_{s}^{t,x,u^{*}},Z_{s}^{t,x,u^{*}})_{t\leq s\leq T}$ is the solution of
BSDE (\ref{stateBSDE}) with optimal feedback control $u^{*}(\cdot)$; }

{ (ii) $(Y_{t}^{t,x,u^{*}},Z_{t}^{t,x,u^{*}})=(g(t,x),\sigma^{*\top}%
\partial_{x}g(t,x))$ is continuous on $[0,T]\times\mathbb{R}^{d}$. }
\end{lemma}

{ Based on Lemma \ref{lem-5}, we introduce a new definition for the viscosity
solution of equations (\ref{hjb-1}) and (\ref{hjb-10}): }

\begin{definition}
{ \label{def-vis-2} Suppose $(w_{1},w_{2})\in C([0,T]\times\mathbb{R}%
^{n})\times C^{0,1}([0,T]\times\mathbb{R}^{d})$. $(w_{1},w_{2})$ is a
viscosity sub-solution of (\ref{hjb-1}) and (\ref{hjb-10}), if $\forall
(\Gamma_{1},\Gamma_{2})\in C^{1,2}([0,T]\times\mathbb{R}^{d})$, for any
$(t,x)\in[0,T]\times\mathbb{R}^{d}$ such that $\Gamma_{i} \geq w_{i}$ and
$\Gamma_{i}(t,x)=w_{i}({t,x}),\ i=1,2$, we have
\begin{equation}%
\begin{array}
[c]{l}%
\displaystyle \underset{u\in{U}}{\inf}\{\mathcal{D}{\Gamma
}_{1}\left(  t,x\right)  +f\left(  t,x,\Gamma_{2}(t,x),\sigma^{\top}%
\partial_{x} \Gamma_{2}(t,x),u\right)  \}\geq0,
\end{array}
\label{subvis-1}%
\end{equation}
and
\begin{equation}%
\begin{array}
[c]{l}%
\displaystyle \mathcal{D}^{*}\Gamma_{2}\left(  t,x\right)  +h\left(
t,x,\Gamma_{2}(t,x),\sigma^{*\top}\partial_{x} \Gamma_{2}(t,x),u^{*}%
(t,x)\right)  \geq0.
\end{array}
\label{subvis-2}%
\end{equation}
$(w_{1},w_{2})$ is a viscosity super-solution of (\ref{hjb-1}) and
(\ref{hjb-10}), if $\forall(\Gamma_{1},\Gamma_{2})\in C^{1,2}([0,T]\times
\mathbb{R}^{d})$, for any $(t,x)\in[0,T]\times\mathbb{R}^{d}$ such that
$\Gamma_{i} \leq w_{i}$ and $\Gamma_{i}(t,x)=w_{i}({t,x}),\ i=1,2$, we have
\begin{equation}%
\begin{array}
[c]{l}%
\displaystyle \underset{u\in{U}}{\inf}\{\mathcal{D}{\Gamma
}_{1}\left(  t,x\right)  +f\left(  t,x,\Gamma_{2}(t,x),\sigma^{\top}%
\partial_{x} \Gamma_{2}(t,x),u\right)  \}\leq0,
\end{array}
\label{supvis-1}%
\end{equation}
and
\begin{equation}%
\begin{array}
[c]{l}%
\displaystyle \mathcal{D}^{*}\Gamma_{2}\left(  t,x\right)  +h\left(
t,x,\Gamma_{2}(t,x),\sigma^{*\top}\partial_{x} \Gamma_{2}(t,x),u^{*}%
(t,x)\right)  \leq0,
\end{array}
\label{supvis-2}%
\end{equation}
$(w_{1},w_{2})$ is a viscosity solution of (\ref{hjb-1}) and (\ref{hjb-10}),
if it is both a viscosity super-solution and a viscosity sub-solution of
(\ref{hjb-1}) and (\ref{hjb-10}) over $[0,T]\times\mathbb{R}^{d}$. }
\end{definition}

{ \label{def-vis-rem}
Clearly, Definition \ref{def-vis-2} of a viscosity solution is stronger than that given in Definition \ref{def-vis}, as it requires the viscosity solution of the auxiliary equation to be differentiable in $x$. However, we do not impose any additional conditions on $w_1$. Therefore, Definition \ref{def-vis-2} reduces to the classical case when the function
$f(t,x,y,z,u)$ is independent of $(y,z)$.}

{ Since $w_{2}\in C^{0,1}([0,T]\times\mathbb{R}^{d})$, $\Gamma_{2}\geq w_{2}$
and $\Gamma_{2}(t,x)=w_{2}({t,x})$, we get that $\partial_{x}\Gamma
_{2}(t,x)=\partial_{x}w_{2}(t,x)$, which is an important observation in the
proof of the following Theorem \ref{the-vis}. This result is also true for the viscosity
super-solution in Definition \ref{def-vis-2}. }

In the following, we present the main result of this section.

\begin{theorem}
{ \label{the-vis} Let Assumptions \ref{ass-1}, \ref{ass-2}, and \ref{ass-4}
hold. Then the value function $(v,g)$ is a viscosity solution of the extended
HJB equation (\ref{hjb-1})--(\ref{hjb-10}). }
\end{theorem}

\noindent\textbf{Proof:} Let $\Gamma_{1},\Gamma_{2} \in C^{1,2}(
[0,T]\times\mathbb{R}^{d})$. For a given $(t,x)\in[0,T]\times\mathbb{R}^{d}$
such that $\Gamma_{1} \geq v,\ \Gamma_{2}\geq g$ and $\Gamma_{1}
(t,x)=v(t,x),\ \Gamma_{2} (t,x)=g(t,x)$, from Remark \ref{def-vis-rem}, we
have $\partial_{x}\Gamma_{2}(t,x)=\partial_{x}g(t,x)$. Now, we prove the
viscosity sub-solution inequalities (\ref{subvis-1}) and (\ref{subvis-2}).

{ \textbf{Step 1:} First we prove that $\Gamma_{2}(t,x)$ satisfies inequality
(\ref{subvis-2}). Applying It\^{o}'s formula to $\Gamma_{2} (s,X^{{t,x},u^{*}}%
_{s})$, for $t< s< T$, we have that
\begin{equation}%
\begin{array}
[c]{cl}%
\mathbb{E}_{t}[\Gamma_{2} (s,X^{{t,x},u}_{s})]-\Gamma_{2}({t,x})
=\mathbb{E}_{t}[\Gamma_{2} (s,X^{{t,x},u^{*}}_{s})]- g({t,x}) =
\displaystyle \mathbb{E}_{t}\int_{t}^{s}\mathcal{D}^{*}\Gamma_{2}
(r,X^{{t,x},u^{*}}_{r})\mathrm{d}r. &
\end{array}
\label{Ito formula0}%
\end{equation}
From Lemma \ref{lem-5}, it follows that
\[
\mathbb{E}_{t}[g(s,X^{t,x,u^{*}}_{s})]-g(t,x)=\mathbb{E}_{t}\left[  \int
_{t}^{s}-h(r,X^{t,x,u^{*}}_{r},Y^{t,x,u^{*}}_{r},Z^{t,x,u^{*}}_{r}%
,u^{*}(r,X^{t,x,u^{*}}_{r}))\mathrm{d}r\right]  .
\]
and $Z^{t,x,u^{*}}_{s}=\sigma^{*\top}g(s,X^{t,x,u^{*}}_{s})$. Combining
equation (\ref{Ito formula0}), we obtain
\[
\displaystyle \mathbb{E}_{t}\int_{t}^{s}\mathcal{D}^{*}\Gamma_{2}
(r,X^{{t,x},u^{*}}_{r})\mathrm{d}r\geq\mathbb{E}_{t}\left[  \int_{t}%
^{s}-h(r,X^{t,x,u^{*}}_{r},Y^{t,x,u^{*}}_{r},Z^{t,x,u^{*}}_{r},u^{*}%
(r,X^{t,x,u^{*}}_{r}))\mathrm{d}r\right]  .
\]
That is,
\[
\displaystyle \mathbb{E}_{t}\int_{t}^{s}\left[  \mathcal{D}^{*}\Gamma_{2}
(r,X^{{t,x},u^{*}}_{r})+h(r,X^{t,x,u^{*}}_{r},Y^{t,x,u^{*}}_{r},Z^{t,x,u^{*}%
}_{r},u^{*}(r,X^{t,x,u^{*}}_{r}))\right]  \mathrm{d}r \geq0.
\]
Letting $s\to t$, it follows that
\[
\mathcal{D}^{*}\Gamma_{2}\left(  t,x\right)  +h\left(  t,x,\Gamma
_{2}(t,x),\sigma^{*\top}\partial_{x} g(t,x),u^{*}(t,x)\right)  \geq0.
\]
Since $\partial_{x}\Gamma_{2}(t,x)=\partial_{x}g(t,x)$, we have
\[
\displaystyle \mathcal{D}^{*}\Gamma_{2}\left(  t,x\right)  +h\left(
t,x,\Gamma_{2}(t,x),\sigma^{*\top}\partial_{x} \Gamma_{2}(t,x),u^{*}%
(t,x)\right)  \geq0.
\]
}

{ \textbf{Step 2:} We now prove that $v$ is the sub-viscosity solution of
(\ref{hjb-1}). Define
\[
\displaystyle u_{1}(r)=%
\begin{cases}
{u(r)}, & \mbox{ $t\leq r< s$},\\
{u}^{\ast}(r,x), & \mbox{ $ s\leq r\leq T$},
\end{cases}
\]
where ${u}^{\ast}(\cdot)$ is an optimal feedback control of value function $v$ on
$(s,T]$. By It\^{o}'s formula, for $t<s<T$, we have that
\begin{equation}%
\begin{array}
[c]{cl}%
\mathbb{E}_{t}[\Gamma_{1}(s,X_{s}^{{t,x},u})]-v({t,x})=\mathbb{E}_{t}%
[\Gamma_{1}(s,X_{s}^{{t,x},u})]-\Gamma_{1}({t,x})=\displaystyle\mathbb{E}%
_{t}\int_{t}^{s}\mathcal{D}\Gamma_{1}(r,X_{r}^{{t,x},u})\mathrm{d}r. &
\end{array}
\label{Ito formula}%
\end{equation}
Based on the DPP established in Theorem \ref{DPP_1}, we have
\[
v({t,x})=\inf_{u(\cdot)\in\mathcal{U}[t,s]}\mathbb{E}_{t}\bigg{[}\int_{t}%
^{s}f(\Theta_{r}^{{t,x},u},u(r))\mathrm{d}r+v(s,X_{s}^{{t,x}%
,u})\bigg{]},
\]
which implies that
\[
v({t,x})\leq\mathbb{E}_{t}\bigg{[}\int_{t}^{s}f(\Theta_{r}^{{t,x},u}%
,u(r))\mathrm{d}r+v(s,X_{s}^{{t,x},u})\bigg{]}.
\]
By equation (\ref{Ito formula}) and observing that for $(t,x)\in
\lbrack0,T]\times\mathbb{R}^{d}$, $\Gamma_{1}\geq v$ and $\Gamma
_{1}(t,x)=v(t,x)$, we deduce that
\[%
\begin{array}
[c]{rl}%
0\leq & \displaystyle\mathbb{E}_{t}\bigg{[}\int_{t}^{s}f(\Theta_{r}^{{t,x}%
,u},u(r))\mathrm{d}r+v(s,X_{s}^{{t,x},u})\bigg{]}-v({t,x})\\
\leq & \displaystyle\mathbb{E}_{t}\bigg{[}\int_{t}^{s}f(\Theta_{r}^{{t,x}%
,u},u(r))\mathrm{d}r+\Gamma_{1}(s,X_{s}^{{t,x},u}%
)\bigg{]}-\Gamma_{1}(t,x)\\
\leq & \displaystyle\mathbb{E}_{t}\int_{t}^{s}\bigg{[}\mathcal{D}\Gamma
_{1}(r,X_{r}^{{t,x},u})+f(\Theta_{r}^{{t,x},u},u(r))\bigg{]}\mathrm{d}r.
\end{array}
\]
Dividing on both sides of the above inequality by $s-t$, we get that
\[%
\begin{array}
[c]{c}%
0\leq\displaystyle\frac{1}{s-t}\displaystyle\mathbb{E}_{t}\int_{t}%
^{s}\bigg{[}\mathcal{D}\Gamma_{1}(r,X_{r}^{{t,x},u})+f(\Theta_{r}^{{t,x}%
,u},u(r))\bigg{]}\mathrm{d}r.
\end{array}
\]
Now, letting $s\downarrow t$, it follows that
\[%
\begin{array}
[c]{c}%
\mathcal{D}{\Gamma}_{1}\left(  t,x\right)  +f\left(  t,x,g(t,x),\sigma^{\top}\partial_{x}g(t,x),u(t)\right)  \geq0.
\end{array}
\]
Since, $\partial_{x}g(t,x)=\partial_{x}\Gamma_{2}(t,x)$, we
have
\[%
\begin{array}
[c]{c}%
\mathcal{D}{\Gamma}_{1}\left(  t,x\right)  +f\left(  t,x,g(t,x),\sigma^{\top}\partial_{x}\Gamma_{2}(t,x),u(t)\right)  \geq0.
\end{array}
\]
Thus
\[
\displaystyle\underset{u\in{U}}{\inf}\{\mathcal{D}{\Gamma
}_{1}\left(  t,x\right)  +f\left(  t,x,\Gamma_{2}(t,x),\sigma^{\top
}\partial_{x}\Gamma_{2}(t,x),u\right)  \}\geq0.
\]
}

\textbf{Step 3:} In the following, we prove the viscosity super-solution. Let
$\Gamma_{1},\Gamma_{2} \in C^{1,2}([0,T]\times\mathbb{R}^{d})$. For given
$(t,x)\in[0,T]\times\mathbb{R}^{d}$ such that $\Gamma_{1} \leq v,\ \Gamma
_{2}\leq g$ and $\Gamma_{1} (t,x)=v(t,x),\ \Gamma_{2} (t,x)=g(t,x)$, from
Remark \ref{def-vis-rem}, we have $\partial_{x}\Gamma_{2} (t,x)=\partial
_{x}g(t,x)$. The proof of $\Gamma_{2}$ satisfying inequality (\ref{supvis-2})
is similar with that in \textbf{Step 1}. Thus, we omit it.

We now use the method of \textquotedblleft proof by
contradiction\textquotedblright to verify that $\Gamma_{1}$ satisfies
inequality (\ref{supvis-1}). Let us assume that, there exist $(t,x)\in
[0,T]\times\mathbb{R}^{d}$, with $\Gamma_{1} \leq v$ on $[0,T]\times
\mathbb{R}^{d}$ and $\Gamma_{1} (t,x)=v(t,x)$, and $\delta_{1}>0$ such that
for any control $u(\cdot)\in\mathcal{U}[t,T]$,
\[
\mathcal{D}{\Gamma}_{1}\left(  t,x\right)  +f\left(  t,x,\Gamma_{2}%
(t,x),\sigma^{\top}\partial_{x} \Gamma_{2}(t,x),u(t)\right)  \geq
2\delta_{1}.
\]

{ Let $u(\cdot)\in\mathcal{U}[t,T]$, and there exists $t<s $ such that
$\mathbb{E}_{t}[\left|  X^{t,x,u}_{s}-x\right|  ] \leq\rho(\delta_{1})$, where
$\rho(\cdot)$ is a continuous function on $\mathbb{R}$ and is independent of
$u$ with $\rho(0)=0$. Thus, for any $r\in\lbrack t,s]$, we have
\begin{equation}%
\begin{array}
[c]{c}%
\mathbb{E}_{t}\bigg{[}\mathcal{D}{\Gamma}_{1}\left(  r,X^{t,x,u}_{r}\right)
+f\left(  r,X^{t,x,u}_{r},g(r,X^{t,x,u}_{r}),\sigma^{*\top}\partial_{x}
\Gamma_{2}(r,X^{t,x,u}_{r}),u(r)\right)  \bigg{]}\geq\delta_{1}.
\end{array}
\label{viscosity-1}%
\end{equation}
}

{ On the other hand, since $\Gamma_{1} \leq v$ on $[0,T]\times\mathbb{R}^{d}$,
by the definition of $J$ and $v$, we have that
\[%
\begin{array}
[c]{rl}%
J(t,x;u_{1}(\cdot))= & \displaystyle \mathbb{E}_{t}\bigg{[}\int_{t}%
^{s}f\left(  r,\Theta^{t,x,u_{1}}_{r},u(r)\right)
\mathrm{d}r+v(s,X^{t,x,u}_{s}) \bigg{]}\\
\geq & \displaystyle \mathbb{E}_{t}\bigg{[}\int_{t}^{s}f\left(  r,\Theta
^{t,x,u_{1}}_{r},u(r)\right)  \mathrm{d}r+\Gamma_{1}%
(s,X^{t,x,u}_{s}) \bigg{]},
\end{array}
\]
where
\[
\displaystyle u_{1}(r)=
\begin{cases}
u(r), & \mbox{ $t\leq r< s$},\\
{u}^{*}(r,x), & \mbox{ $ s\leq r\leq T$}.
\end{cases}
\]
It follows equation (\ref{viscosity-1}) that
\begin{equation}%
\begin{array}
[c]{cl}%
J(t,x;u_{1}(\cdot))\geq & \displaystyle\mathbb{E}_{t}\int_{t}^{s}%
\bigg{[}\delta_{1}-\mathcal{D}{\Gamma}_{1}\left(  r,X^{t,x,u}_{r}\right)
\bigg{]}ds+\mathbb{E}_{t}\bigg{[}\Gamma_{1} (s,X^{{t,x},u}_{s})\bigg{]}.
\end{array}
\label{viscosity-2}%
\end{equation}
In addition, similar to (\ref{Ito formula}), it follows that
\begin{equation}%
\begin{array}
[c]{cl}%
\mathbb{E}_{t}[\Gamma_{1} (s,X^{{t,x},u}_{s})]-v({t,x}) =\mathbb{E}_{t}%
[\Gamma_{1} (s,X^{{t,x},u}_{s})]-\Gamma_{1} ({t,x}) = \displaystyle \mathbb{E}%
_{t}\int_{t}^{s}\mathcal{D}\Gamma_{1} (r,X^{{t,x},u}_{r})\mathrm{d}r. &
\end{array}
\label{viscosity-3}%
\end{equation}
Therefore, from equations (\ref{viscosity-2}) and (\ref{viscosity-3}), we
have
\[
J(t,x;u_{1}(\cdot))\geq(s-t)\delta_{1}+\Gamma_{1} ({t,x})=(s-t)\delta
_{1}+v({t,x}).
\]
Taking the minimum over $u(t,x)\in\mathcal{U}_{t}$, we have that
\[
v ({t,x})\geq(s-t)\delta_{1}+v ({t,x}).
\]
This contradicts to the fact that $\delta_{1}>0$. Therefore, $v({t,x})$ is a
viscosity super-solution of equation (\ref{hjb-1}). $\ \ \ \ \ \ \ \Box$ }

\begin{theorem}
{ \label{the-vis-un} Let Assumptions \ref{ass-1}, \ref{ass-2} and \ref{ass-4}
hold. Then the value function $(v,g)$ is the unique viscosity solution of the
HJB equations (\ref{hjb-1}) and (\ref{hjb-10}). Namely, }

{ (i) For a given optimal feedback control $u^{*}(\cdot)$ of value function
(\ref{value-v}), the function (\ref{value-g}) is the unique viscosity solution
of (\ref{hjb-10}); }

{ (ii) The value function (\ref{value-v}) is the unique viscosity solution of
(\ref{hjb-1}). }
\end{theorem}

{ \noindent\textbf{Proof:} For a given optimal feedback control $u^{\ast}(\cdot)$ of
value function (\ref{value-v}), we first consider equation (\ref{hjb-10}) with
$g(T,x)=\Phi\left(  x\right)  $. Note that the new Definition \ref{def-vis-2}
of viscosity solution is stronger than Definition \ref{def-vis}. Following
the classical Definition \ref{def-vis} of viscosity solution, we can obtain
the uniqueness of viscosity solution of (\ref{hjb-10}) based on the results of
\citet{CIL92}, which implies the uniqueness of viscosity solution in
Definition \ref{def-vis-2}. }

{ Now, from Lemma \ref{lem-5} and Theorem \ref{the-vis}, the value function
$g(t,x)\in C^{0,1}([0,T]\times\mathbb{R}^{d})$ of (\ref{value-g}) is the
unique viscosity solution of (\ref{hjb-10}) under a given optimal control
$u^{*}(\cdot)$. {We now consider equation (\ref{hjb-1}) for any given}
$(t,x)\in[0,T]\times\mathbb{R}^{d}$. Note that $g(t,x)\in C^{0,1}%
([0,T]\times\mathbb{R}^{d})$. Hence the second item in equation (\ref{hjb-1})
are continuous in $(t,x)$. Therefore, the uniqueness of the viscosity solution
in Definition \ref{def-vis-2} is obtained directly by \citet{CIL92}. Then,
from Theorem \ref{the-vis}, the value function $v(t,x)$ of (\ref{value-v}) is
the unique viscosity solution of (\ref{hjb-1}). }

{ For another optimal feedback control $u^{\ast}(\cdot)$ and the related $g(t,x)$,
repeat the above procedure, we can still prove that the value function
$v(t,x)$ of (\ref{value-v}) is the unique viscosity solution of (\ref{hjb-1}).
Hence the uniqueness of viscosity solution to (\ref{value-v}) does not depend
on either the choice of an optimal control or the related $g(t,x)$.
$\ \ \ \ \ \ \ \Box$ }

\section{Connections With Time-Inconsistent Control Problems}
\label{sec:time-in}

In this section, we first demonstrate a relation between problem (\ref{in-value-3}) and
\citet{P93}. The relation of our results with the equilibrium concept in
\citet{BKM2017} is also investigated. 

\subsection{\bigskip Relationship with \citet{P93}}

\citet{P93} investigated a general control problem with cost functional:%
\begin{equation}
{\tilde{J}}(t,x;u(\cdot))=\mathbb{E}_{t}\left[  \int_{t}^{T}f\left(
s,X_{s}^{t,x,u},Y_{s}^{t,x,u},Z_{s}^{t,x,u},u(s)\right)
\mathrm{d}s+\gamma(Y_{t}^{t,x,u})+G(X_{T}^{t,x,u})\right]  , \label{6cost-0}%
\end{equation}
where $\left(  X^{t,x,u}\right)  $ and $\left(  Y^{t,x,u},Z^{t,x,u}\right)  $
satisfy SDE \eqref{state-1} and BSDE \eqref{state2} respectively. A nonlinear
term $\gamma(Y_{t}^{t,x,u})$ was first introduced in the cost functional. Note
that the control problem%
\begin{equation}
\tilde{v}(t,x)=\inf_{u\in\mathcal{U}[t,T]}{\tilde{J}}(t,x;u(\cdot))
\label{6value-v}%
\end{equation}
may be time inconsistent due to the appearance of $\gamma(Y_{t}^{t,x,u})$.
\citet{P93} established a local stochastic maximum principle at initial time
$t=0$ for a controlled FBSDE system. However, an HJB characterization to
problem \eqref{6value-v} is still open until now. The difficulty
lies in the possible time-inconsistency. Interestingly, we find that, if the
function $\gamma(\cdot)$ is twice differentiable, problem \eqref{6value-v} can
be converted into a multidimensional time-consistent one.

In fact, it follows the cost functional (\ref{6cost-0}) that
\begin{equation}
{\tilde{J}}(s,X_{s}^{t,x,u};u(\cdot))=\mathbb{E}_{s}\left[  \int_{s}%
^{T}f\left(  r,X_{r}^{t,x,u},Y_{r}^{t,x,u},Z_{r}^{t,x,u},u(r)\right)  \mathrm{d}s+\gamma(Y_{s}^{t,x,u})+G(X_{T}^{t,x,u})\right]  .
\label{6cost-01}%
\end{equation}
Combining (\ref{6cost-0}) and (\ref{6cost-01}), we have
\begin{equation}%
\begin{array}
[c]{rl}%
{\tilde{J}}(t,x;u(\cdot))= & \mathbb{E}_{t}\big[\int_{t}^{s}f\left(
r,X_{r}^{t,x,u},Y_{r}^{t,x,u},Z_{r}^{t,x,u},u(r)\right)
\mathrm{d}r\\
& +\gamma(Y_{t}^{t,x,u})-\gamma(Y_{s}^{t,x,u})+{\tilde{J}}(s,X_{s}%
^{t,x,u};u(\cdot))\big].
\end{array}
\label{6cost-02}%
\end{equation}
Assuming $\gamma(\cdot)$ has continuous second-order derivatives on
$\mathbb{R}$ and applying It\^{o}'s formula to $\gamma(Y_{t}^{t,x,u})$, we get
that,
\[
\gamma(Y_{t}^{t,x,u})=\mathbb{E}_{t}\left[  \int_{t}^{T}\left(  h(s)\gamma
_{y}(Y_{s}^{t,x,u})-\frac{1}{2}Z_{s}^{t,x,u}({Z_{s}^{t,x,u}})^{\top}%
\gamma_{yy}(Y_{s}^{t,x,u})\right)  \mathrm{d}s+\gamma(\Phi(X_{T}%
^{t,x,u}))\right]  ,
\]
and
\begin{equation}
\mathbb{E}_{t}\left[  \gamma(Y_{t}^{t,x,u})-\gamma(Y_{s}^{t,x,u})\right]
=\mathbb{E}_{t}\left[  \int_{t}^{s}\left(  h(r)\gamma_{y}(Y_{r}^{t,x,u}%
)-\frac{1}{2}Z_{r}^{t,x,u}({Z_{r}^{t,x,u}})^{\top}\gamma_{yy}(Y_{r}%
^{t,x,u})\right)  \mathrm{d}r\right]  , \label{6adjust}%
\end{equation}
where $h(r):=h(r,X_{r}^{t,x,u},Y_{r}^{t,x,u},Z_{r}^{t,x,u},u(r,X_{r}%
^{t,x,u}))$. The right hand of \eqref{6adjust} can be regarded as an adjustment between $\gamma
(Y_{t}^{t,x,u})$ and $\gamma(Y_{s}^{t,x,u})$ from time $t$ to time $s$, such
that there exists a dynamic time-consistent structure from $\gamma(Y_{t}^{t,x,u})$ to
$\gamma(Y_{s}^{t,x,u})$. {It is important to note that the transformation via  It\^{o}'s formula is applied at a fixed time point. However, the transformed problem acquires a new dynamic structure, which may potentially alter the original time-inconsistency property.}

Denote
\[
\tilde{f}(t,x,y,z,u)=f(t,x,y,z,u)+h(t)\gamma_{y}(y)-\frac{1}{2}zz^{\top}%
\gamma_{yy}(y),\ \ \text{and}\quad\tilde{G}(x)=G(x)+\gamma(\Phi(x)).
\]
It follows \eqref{6cost-01} and \eqref{6cost-02} that
\begin{align}
{\tilde{J}}(t,x;u(\cdot))  &  =\mathbb{E}_{t}\left[  \int_{t}^{T}\tilde
{f}\left(  s,X_{s}^{t,x,u},Y_{s}^{t,x,u},Z_{s}^{t,x,u},u(s)\right)  \mathrm{d}s+\tilde{G}(X_{T}^{t,x,u})\right] \label{tc6.1}\\
&  =\mathbb{E}_{t}\left[  \int_{t}^{s}\tilde{f}\left(  r,X_{r}^{t,x,u}%
,Y_{r}^{t,x,u},Z_{r}^{t,x,u},u(r)\right)  \mathrm{d}r+\tilde
{J}(s,X_{s}^{t,x,u};u(\cdot))\right]  .\nonumber
\end{align}
{It is precisely this  structural adjustment that ensures the original problem becomes time-consistent under the reformulated running cost structure.}

Then applying Theorem \ref{DPP_1}, we have the DPP:
\[
\tilde{v}({t,x})=\displaystyle\inf_{u(\cdot)\in\mathcal{U}[t,s]}\mathbb{E}%
_{t}\bigg{[}\int_{t}^{s}\tilde{f}(r,X_{r}^{t,x,u},Y_{r}^{t,x,u},Z_{r}%
^{t,x,u})\mathrm{d}r+\tilde{v}(s,X_{s}^{{t,x},{u}})\bigg{]}.
\]
The above results are concluded as follows:

\begin{proposition}
\label{6pro-1} Let Assumptions \ref{ass-1} and \ref{ass-2} hold, and
$\gamma(\cdot)$ be twice differentiable on $\mathbb{R}$. Then the extended DPP
holds for the optimality problem \eqref{tc6.1} and ${\tilde{v}}(t,x)$ satisfies the
extended HJB equation \eqref{hjb-1} and \eqref{hjb-10} with $\left(
f,G\right)  $ replaced by ($\tilde{f},\tilde{G}).$

\end{proposition}
{Therefore, Proposition \ref{6pro-1}  provides a PDE characterization for the general class of control problems considered in \citet{P93}, thereby filling the long-standing gap in the literature regarding the absence of such a characterization.}
In the next subsection \ref{subsec:comparison}, we will show that our extended DPP provides an equilibrium solution to the original problem \eqref{6value-v}, in the sense of  \textquotedblleft Nash subgame perfect equilibrium\textquotedblright\  (\citet{BM2010}).

\subsection{Comparison to {\citet{BKM2017}}}
\label{subsec:comparison}
In the draft version \citet{BM2010} and the published paper {\citet{BKM2017},}
they investigated the optimization of (\ref{6cost-0}) with $Y^{t,x,u}$ being a
conditional expectation and $f$ not depending on $Z^{t,x,u}$. They showed that
the optimization problem for maximizing the reward function (\ref{6cost-0})
does not satisfy the Bellman optimality principle due to the dependence on a
nonlinear function of the expectation. Thus it is a time-inconsistent problem.
Classical dynamic programming is therefore not available for solving this
problem. To give a time-consistent solution, they formulated this
time-inconsistent problem in a game theoretic framework. The equilibrium, if
exists, is what we want to look for the ``optimal control". More precisely,
the optimization problem is treated as a non-cooperative game and at each
point of time $t$, there is a player $t$ which can be regarded as an
incarnation of the investor. Then the equilibrium strategy
${{u}^{\ast}}$ is defined as: for an arbitrary time $t$, the optimal strategy
for player $t$ is ${{u}}^{\ast}(t,\cdot)$, supposing that each player $s$
where $s>t$ uses the strategy ${{u}^{\ast}}(s,\cdot)$. The following
mathematical characterization of an equilibrium comes from {\citet{BKM2017}}.

\begin{definition}
\label{def-equi} We call ${{u}^{\ast}}$ an equilibrium strategy of value function
(\ref{6value-v}), if for every admissible control ${u}$ valued in
$\mathbb{R}^{d}$ and $h>0$,
\[
{u}_{h}(s,x)=%
\begin{cases}
{u}, & \mbox{for \ensuremath{t\le s<t+h,\quad x\in \mathbb{R}^n}}\\
{{u}^{\ast}(s,x)}, & \mbox{for \ensuremath{t+h\le s\le T,\quad x\in
\mathbb{R}^n}},
\end{cases}
\]
such that
\begin{equation}
\limsup_{h\rightarrow0^{+}}\frac{\tilde{J}(t,x;{{u}^{\ast}(\cdot)})-\tilde
{J}(t,x;{u}_{h}(\cdot))}{h}\leq0 \label{equ-1}%
\end{equation}
for any $(t,x)\in\lbrack0,T]\times\mathbb{R}^{d}$.
\end{definition}

From the above definition, when $h$ is small enough, we have $\tilde
{J}(t,x;{u}^{\ast}(\cdot))-\tilde{J}(t,x;{u}_{h}(\cdot))\leq0$, which is
consistent with the notion of Nash equilibrium for non-cooperative games with finite number of players. Now
we show that the time-consistent control provided by Proposition {\ref{6pro-1}
}is an equilibrium.

\begin{proposition}
\label{the6-0} Let Assumptions \ref{ass-1} and \ref{ass-2} hold. Then, for any
given $(t,x)\in\lbrack0,T)\times\mathbb{R}^{d}$, the dynamic consistent
optimal control $u^{*}(\cdot)$ of the optimality problem \eqref{tc6.1} given in
Proposition \ref{6pro-1} is an equilibrium strategy of the original problem (\ref{6value-v}).
\end{proposition}

{ \noindent\textbf{Proof:} Let $u^{\ast}(\cdot)$ be an time-consistent optimal
control of the optimality problem \eqref{tc6.1}. By Proposition \ref{6pro-1}, we
have that
\begin{align}
\tilde{v}({t,x})&  =\displaystyle\inf_{u(\cdot)\in\mathcal{U}[t,s]}\mathbb{E}%
_{t}\bigg{[}\int_{t}^{s}\tilde{f}(r,X_{r}^{t,x,u},Y_{r}^{t,x,u},Z_{r}%
^{t,x,u})\mathrm{d}r+\tilde{v}(s,X_{s}^{{t,x},{u}})\bigg{]}\nonumber\\
&  =\tilde{J}(t,x;u^{\ast}(\cdot))\leq\tilde{J}(t,x;u_{h}(\cdot)).\nonumber
\end{align}
Consequently,
\[
\limsup_{h\rightarrow0^{+}}\frac{\tilde{J}(t,x;{{u}^{\ast}(\cdot)})-\tilde
{J}(t,x;{u}_{h}(\cdot))}{h}\leq0.
\]
Thus $u^{\ast}(\cdot)$ is an equilibrium strategy of value function
(\ref{6value-v}). $\ \ \ \ \ \ \ \Box$ }

\bigskip An important contribution of \citet{BKM2017} is to extend the dynamic
mean-variance model in \citet{Basak2010} to a fair general framework. However,
they did not provide a DPP. The HJB system is informally derived and a
verification theorem is given when the HJB system has a smooth solution. The
uniqueness and existence of viscosity solution for the extended HJB is also
left open, which is regarded as ``technically extremely difficult".
See\ conclusions in their paper \citet{BKM2017} for these open problems. Our
paper provides some answers from the following aspects. First, we transform
the original problem into a multidimensional time-consistent control problem.
Thus a DPP is available to derive the corresponding HJB system. The above
Proposition \ref{the6-0} shows that our time-consistent optimal control is
just an equilibrium of the original problem. Furthermore, viscosity solution
is also studied in this paper. The uniqueness and existence of viscosity solutions
for the extended HJB are investigated in both cases, whether $f$ depends on $Z^{t,x,u}$ or not.
Note that the uniqueness of a solution to a multidimensional PDE depends heavily on the comparison theorem.
Thus, it is usually challenging to identify the uniqueness of a solution to a multidimensional PDE. In conclusion, our results not only address some open problems proposed in \citet{BKM2017}, but also extend their model to the extent where the conditional expectation can be replaced by a nonlinear expectation, i.e., a solution to a BSDE, particularly, $f$ can also depend on $Z^{t,x,u}$.

\section{Applications and Examples}

\label{sec:time-exam} In this section, we provide three
examples to verify our main results, including the classical mean-variance
model, utility optimization for a narrow framing investor and risk-sensitive
portfolio optimization under a nonlinear expectation.

\subsection{The mean-variance model}

\label{mean-variance}We start with a market in which two assets are available
for investment: a riskless asset with constant interest rate $r>0$, and
a risky asset. The stock price evolves according to the following
geometric Brownian motion:
\[
\mathrm{d}S_{t}=\mu S_{t}\mathrm{d}t+\sigma S_{t}\mathrm{d}B_{t},
\]
where $B_{t}$ is a standard 1-dimensional Brownian motion, the constants $\mu
$, $\sigma$ represent the expected return rate and the volatility
respectively.

A self-financing wealth process $\left(  X_{s}\right)  _{s\geq t}$ can be described by
\begin{equation}
\mathrm{d}X_{s}^{t,x,u}=\left[  rX_{s}^{t,x,u}+\left(  \mu-r\right)  u(s)\right]  \mathrm{d}s+\sigma u(s)\mathrm{d}B_{s}, \label{wealth}%
\end{equation}
where $X_{t}^{t,x,u}=x$ and $u(\cdot)$ is the dollar amount invested in the
stock, standing for a trading strategy.

In the traditional mean-variance model, we take the reward function as
\[
{\tilde{J}}(t,x;u(\cdot))=\mathbb{E}_{t}[X_{T}^{t,x,u}]-\frac{\lambda}%
{2}\mathbb{E}_{t}[(X_{T}^{t,x,u}-\mathbb{E}_{t}[X_{T}^{t,x,u}])^{2}],
\]
which is the well-known dynamic mean-variance problem (\citet{Basak2010},
\citet{BKM2017}), where $\lambda>0$ is the risk aversion
coefficient. We rewrite ${\tilde{J}}(t,x;u(\cdot))$ as follows
\[
{\tilde{J}}(t,x;u(\cdot))=\mathbb{E}_{t}[X_{T}^{t,x,u}-\frac{\lambda}{2}%
(X_{T}^{t,x,u})^{2}]+\frac{\lambda}{2}\left(  Y_{t}^{t,x,u}\right)  ^{2},
\]
where
\begin{equation}
Y_{t}^{t,x,u}=X_{T}^{t,x,u}-\int_{t}^{T}Z_{s}^{t,x,u}\mathrm{d}B_{s}%
.\label{cond_exp}%
\end{equation}

By the method developed in Section \ref{sec:time-in}, we obtain that
\[
{\tilde{J}}(t,x;u(\cdot))=\mathbb{E}_{t}\left[  X_{T}^{t,x,u}-\frac{\lambda
}{2}\int_{t}^{T}({Z_{s}^{t,x,u}})^{2}\mathrm{d}s\right]  ,
\]
and the value function $\tilde{v}(t,x)$ satisfies the following DPP:
\begin{equation}
\tilde{v}(t,x)=\displaystyle\sup_{u\in\mathcal{U}[t,s]}\mathbb{E}_{t}\left[
\tilde{v}(s,X_{s}^{t,x,u})-\frac{\lambda}{2}\int_{t}^{s}({Z_{r}^{t,x,u}}%
)^{2}\mathrm{d}r\right]  .\label{dpp1-7.1}%
\end{equation}

\citet{Basak2010} showed that the value function $\tilde{v}(t,x)$ satisfies
the following recursive equation
\begin{equation}
\tilde{v}(t,x)=\displaystyle\sup_{u\in\mathcal{U}[t,s]}\mathbb{E}_{t}\left[
\tilde{v}(s,X_{s}^{t,x,u})-\frac{\lambda}{2}\left(  \mathbb{E}_{t+s}%
[X_{T}^{t,x,u}]-\mathbb{E}_{t}[X_{T}^{t,x,u}]\right)  ^{2}\right]
,\label{dpp2-7.1}%
\end{equation}
which is equivalent to \eqref{dpp1-7.1}. In fact, by \eqref{cond_exp} we
deduce that%
\[
\mathbb{E}_{t}\left(  \mathbb{E}_{t+s}[X_{T}^{t,x,u}]-\mathbb{E}_{t}%
[X_{T}^{t,x,u}]\right)  ^{2}=\mathbb{E}_{t}\left(  Y_{t+s}^{t,x,u}%
-Y_{t}^{t,x,u}\right)  ^{2}=\mathbb{E}_{t}\int_{t}^{s}({Z_{r}^{t,x,u}}%
)^{2}\mathrm{d}r.
\]
Thus, the item $\mathbb{E}_{t}\int_{t}^{s}({Z_{r}^{t,x,u}})^{2}\mathrm{d}r$ is
an adjustment caused by variance. Therefore, our DPP is in line with the one
in \citet{Basak2010} for the mean-variance problem, and generalizes it to a
more general framework.

Furthermore, by Theorem \ref{the-class}, the value function $(\tilde{v}(t,x),g(t,x))$ satisfies the
following extended HJB equation%
\begin{equation}%
\begin{cases}
\underset{u\in{U}}{\sup}\{\mathcal{D}{\tilde{v}}\left(
t,x\right)  -\frac{\lambda}{2}\left(  \partial_{x}g\left(  t,x\right)  \sigma
u\right)  ^{2}\}=0, & \\
\text{ \ \ \ }\mathcal{D}^{\ast}g\left(  t,x\right)  =0, &
\end{cases}
\label{hjb_mv}%
\end{equation}
subject to ${\tilde{v}}(T,x)=x$ and $g(T,x)=x$, where $g(t,x):=Y_{t}%
^{t,x,u^{\ast}}$.
The solution of optimal control $u^{\ast}(\cdot)$ and value function are given as%
\begin{equation}%
\begin{cases}
u^{\ast}(t,x)=\frac{\mu-r}{\lambda\sigma^{2}}e^{-r(T-t)}, & \\
\tilde{v}(t,x)=xe^{r(T-t)}+\frac{(\mu-r)^{2}}{2\lambda\sigma^{2}}(T-t), &
\end{cases}
\label{control_mv}%
\end{equation}
which coincides with the results given in \citet{Basak2010} and
{\citet{BKM2017}}. 

\subsection{Utility optimization for a narrow framing investor}

Suppose $X^{t,x,u}$ is a portfolio following equation \eqref{wealth} which can
be hedged by a risk-free asset and a risky asset. In finance or insurance, if
$X_{T}^{t,x,u}$ is the terminal surplus, then $\mathbb{E}X_{T}^{t,x,u}$
usually represents the cost (price) of the issued product. Consider a narrow
framing investor who would like to maximize local utility on the difference of
surplus deducting part of issued cost, namely, to maximize the following
utility functional:
\[
J(t,x;u(\cdot))=\mathbb{E}_{t}[U(X_{T}^{t,x,u}-\alpha\mathbb{E}_{t}%
X_{T}^{t,x,u})],
\]
where $\alpha\in\left[  0,1\right]  $ is the discount rate. When $\alpha=0$,
the investor maximizes utility on wealth; when $\alpha=1$, the investor plays
a gamble on the discount-free price. Narrow framing means that an individual
uses a utility function that depends directly on the outcome of the gamble
when deciding whether to accept a gamble. See \citet{KL1993}, \citet{K2003},
\citet{RW2009}, \citet{CZZ2022} and reference therein for the notion of narrow framing.

Consider the CARA utility: $U(x)=-{\frac{1}{\lambda}}e^{-\lambda x}$, where
$\lambda>0$. Let $\left(  Y_{1},Z_{1}\right)  $ and $\left(  Y_{2}%
,Z_{2}\right)  $ be solutions of the following two BSDEs, respectively:%
\begin{align}
Y_{1s}^{t,x,u} &  =\alpha X_{T}^{t,x,u}-\int_{s}^{T}Z_{1r}^{t,x,u}%
\mathrm{d}B_{r}\text{, }\\
Y_{2s}^{t,x,u} &  =-e^{-\lambda X_{T}^{t,x,u}}-\int_{s}^{T}Z_{2r}%
^{t,x,u}\mathrm{d}B_{r}\text{.\ \ }%
\end{align}
Then the reward function is
\begin{align}
J(t,x;u(\cdot)) &  =\mathbb{E}_{t}[U(X_{T}^{t,x,u}-\alpha\mathbb{E}_{t}%
X_{T}^{t,x,u})]=-U\left(  -Y_{1t}^{t,x,u}\right)  \cdot Y_{2t}^{t,x,u}%
\nonumber\\
&  =\mathbb{E}_{t}\left[  \int_{t}^{T}-\left(  e^{\lambda Y_{1s}^{t,x,u}%
}Z_{1s}^{t,x,u}Z_{2s}^{t,x,u}+\frac{\lambda e^{\lambda Y_{1s}^{t,x,u}}}%
{2}Y_{2s}^{t,x,u}(Z_{1s}^{t,x,u})^{2}\right)  \mathrm{d}s\right]  -\frac
{1}{\lambda}.
\end{align}
Define $v(t,x):=\max_{u(\cdot)\in\mathcal{U}[t,T]}J(t,x;u(\cdot))$ and
$g_{1}(t,x):=$ $Y_{1t}^{t,x,u^{\ast}}$, $g_{2}(t,x):=$ $Y_{2t}^{t,x,u^{\ast}}%
$. Then by Theorem \ref{the-class}, we get the following extended HJB
equation:%
\begin{equation}%
\begin{cases}
\underset{u\in{U}}{\sup}\{\mathcal{D}{v}\left(  t,x\right)
-\left(  \sigma u\right)  ^{2}\left[  e^{\lambda g_{1}}\cdot\partial_{x}%
g_{1}\cdot\partial_{x}g_{2}+\frac{\lambda e^{\lambda g_{1}}}{2}g_{2}%
\cdot(\partial_{x}g_{1})^{2}\right]  \}=0, & v(T,x)=-{\frac{1}{\lambda}}\\
\text{ \ \ \ \ \ }\mathcal{D}^{\ast}g_{1}\left(  t,x\right)  =0, &
g_{1}(T,x)=\alpha x\\
\text{ \ \ \ \ \ }\mathcal{D}^{\ast}g_{2}\left(  t,x\right)  =0, &
g_{2}(T,x)=-e^{-\lambda x}%
\end{cases}
\label{pde_u1}%
\end{equation}
The infinitesimal operator $\mathcal{D}$ is $\mathcal{D}(\cdot)=[\partial
_{t}+\left(  rx+\left(  \mu-r\right)  u\right)  \partial_{x}+\frac{1}%
{2}\left(  \sigma u\right)  ^{2}\partial_{xx}](\cdot)$, $\mathcal{D}^{\ast
}(\cdot)=[\partial_{t}+\left(  rx+\left(  \mu-r\right)  u^{\ast}\right)
\partial_{x}+\frac{1}{2}\left(  \sigma u^{\ast}\right)  ^{2}\partial
_{xx}](\cdot)$.

By the first-order condition, the optimal control is
\[
u^{\ast}\left(  t,x\right)  =\frac{\mu-r}{\sigma^{2}}\cdot\frac{\partial_{x}%
v}{-\partial_{xx}v+2e^{\lambda g_{1}}\cdot\partial_{x}g_{1}\cdot\partial
_{x}g_{2}+\lambda{e^{\lambda g_{1}}}\cdot g_{2}\cdot(\partial_{x}g_{1})^{2}}.
\]
Assume
\begin{align*}
g_{1}\left(  t,x\right)   &  =a\left(  t\right)  x+b\left(  t\right)  \text{,
}\\
g_{2}\left(  t,x\right)   &  =c\left(  t\right)  e^{d\left(  t\right)
x}\text{,}\\
v\left(  t,x\right)   &  =A\left(  t\right)  e^{B\left(  t\right)
x}\text{.\ }%
\end{align*}
Substituting the above structures into $u^{\ast}$, we get that
\[
u^{\ast}\left(  t\right)  =\frac{\mu-r}{\sigma^{2}}\cdot\frac{AB}%
{-AB^{2}+2acde^{\lambda b}+\lambda{a}^{2}c{e^{\lambda b}}},
\]
which is a function of time $t$. Then substituting $u^{\ast}\left(  t\right)
$ into equation \eqref{pde_u1}, we obtain the following ODE system:%
\begin{equation}%
\begin{cases}
a'+ra=0, & a(T)=\alpha\\
b'+\left(  \mu-r\right)  au^{\ast}=0, & b(T)=0\\
c'+\left(  \left(  \mu-r\right)  du^{\ast}\left(  t\right)  +\frac
{\sigma^{2}}{2}d^{2}\left(  u^{\ast}\right)  ^{2}\right)  c=0, & c(T)=-1\\
d'+ra=0, & d(T)=-\lambda\\
A'+\left(  \mu-r\right)  u^{\ast}BA+\frac{\sigma^{2}}{2}\left(  u^{\ast
}\right)  ^{2}\left(  AB^{2}-2acde^{\lambda b}-{\gamma a}^{2}c{e^{\lambda b}%
}\right)  =0, & A(T)=-\frac{1}{\lambda}\\
B'+rB=0, & B(T)=\left(  \alpha-1\right)  \lambda
\end{cases}
\end{equation}

By a quite tedious calculation, we solve that
\begin{align*}
a\left(  t\right)   &  =\lambda e^{r\left(  T-t\right)  }\text{, }\\
b\left(  t\right)   &  =\frac{\alpha\left(  1-\alpha\right)  \theta^{2}%
}{\lambda}\left(  T-t\right)  ,\\
c\left(  t\right)   &  =-e^{-\frac{\left(  1-\alpha^{2}\right)  \theta^{2}}%
{2}\left(  T-t\right)  }\text{,}\\
d\left(  t\right)   &  =-\lambda e^{r\left(  T-t\right)  }\text{,}\\
A\left(  t\right)   &  =-\frac{1}{\lambda}e^{-\frac{\left(  1-\alpha\right)
^{2}\theta^{2}}{2}\left(  T-t\right)  }\text{,}\\
B\left(  t\right)   &  =\left(  \alpha-1\right)  \lambda e^{r\left(
T-t\right)  }\text{,}%
\end{align*}
where we denote $\theta=\frac{\mu-r}{\sigma}$. Therefore, the optimal feedback control
is
\begin{equation}
u^{\ast}\left(  t\right)  =\frac{\mu-r}{{\lambda}\sigma^{2}}e^{-r\left(
T-t\right)  }\left(  1-\alpha\right)  .\label{control_u1}%
\end{equation}
When $\alpha=0$, the investor only maximizes utility on wealth and his optimal
control is the classical one $u^{\ast}\left(  t\right)  =\frac{\mu-r}%
{{\lambda}\sigma^{2}}e^{-r\left(  T-t\right)  }$; when $\alpha$ increases, the
investor reduces his risky investment to zero until $\alpha=1$, which implies
that one would not participate a risky investment if his reference is the average
level of wealth.

\subsection{Risk-sensitive portfolio optimization under nonlinear expectation}

Risk-sensitive optimal control problem can be traced back to \citet{BS85}.
Then many applications and extensions are developed, see, for instance,
\citet{W90} for an application to differential games and \citet{J73} for a
comprehensive book. Financial application of risk-sensitive control was
pioneered by \citet{LM94} and \citet{BP99}. In particular, \citet{BP99}
proposed the logarithm of the investor's wealth as a reward function, so that
the investor's objective is to maximize the risk-sensitive (log) return of
his/her portfolio or alternatively to maximize a function of the power utility
of terminal wealth. \citet{Hansen2006} connected the risk-sensitive objective
to a robust criteria in which perturbations are characterized by the relative
entropy. However, what they considered are the precommitted strategy at a fixed time.

We consider a risk-sensitive portfolio optimization problem with the reward
functional defined by
\begin{equation}
J(t,x;u(\cdot))=-\frac{1}{\lambda}\ln\mathbb{E}_{t}^{g}[\exp\left(  -\lambda
X_{T}^{t,x,u}\right)  ],\label{risksens}%
\end{equation}
where $\lambda$ is interpreted as a risk-sensitivity parameter or risk
aversion parameter, and
\[
\mathrm{d}X_{s}^{t,x,u}=\mu(s,X_{s}^{t,x,u},u(s))\mathrm{d}%
s+\sigma(s,X_{s}^{t,x,u},u(s))\mathrm{d}B_{s}\text{, \ \ }%
X_{t}^{t,x,u}=x\text{,}%
\]
and $g$-expectation (\cite{P97}) $\mathbb{E}_{t}^{g}[\cdot]:=Y_{t}^{t,x,u}$ is the solution
of the following BSDE at time $t$,
\[
\mathrm{d}Y_{s}^{t,x,u}=-h(X_{s}^{t,x,u},Y_{s}^{t,x,u},Z_{s}^{t,x,u}%
)\mathrm{d}s+Z_{s}^{t,x,u}\mathrm{d}B_{s},\quad Y_{T}^{t,x,u}=\exp\left(
-\mathcal{\lambda}X_{T}^{t,x,u}\right)  .
\]
Consequently, we can rewrite the risk-sensitive reward functional as
\[
J(t,x;u(\cdot))=-\frac{1}{\lambda}\ln Y_{t}^{t,x,u}.
\]
Note that the risk-sensitive portfolio optimization is time-inconsistent since we don't have a DPP with log-expectation.
The usual approach is to apply an exponential transformation to the reward function, transforming it into a classical control problem. This allows for the study of pre-commitment strategies.

Our method can provide equilibrium solutions for risk-sensitive problems.
Applying It\^{o}'s formula to $-\frac{1}{\lambda}\ln Y_{t}^{t,x,u}$, we have
\begin{equation}
J(t,x;u(\cdot))=X_{T}^{t,x,u}+\int_{t}^{T}\left[  \frac{-h_{s}}{\lambda
Y_{s}^{t,x,u}}-\frac{(Z_{s}^{t,x,u})^{2}}{2\lambda(Y_{s}^{t,x,u})^{2}}\right]
\mathrm{d}s+\int_{t}^{T}\frac{Z_{s}^{t,x,u}}{\lambda Y_{s}^{t,x,u}}%
\mathrm{d}B_{s},\label{514}
\end{equation}
where $h_{s}:=h(X_{s}^{t,x,u},Y_{s}^{t,x,u},Z_{s}^{t,x,u})$. Therefore, the
time-consistent structure of (\ref{514}) is
a particular case of our model (\ref{in-value-3}). Let $\tilde{v}%
(t,x)=\sup_{u\in\mathcal{U}[t,T]}J(t,x;u(\cdot))$, and $g(t,x)=Y_{t}%
^{t,x,u^{\ast}}$. Thus the extended HJB equation is%
\begin{equation}%
\begin{cases}
\underset{u\in{U}}{\sup}\{\mathcal{D}\tilde{v}\left(
t,x\right)  +\frac{h(t,x,g(t,x),\partial_{x}g(t,x)\sigma(t,x,u%
))}{\lambda g(t,x)}-\frac{(\partial_{x}g(t,x)\sigma(t,x,u))^{2}}{2\lambda
g^{2}(t,x)}\}=0, & \\
\text{ \ \ \ \ }\mathcal{D}^{\ast}g\left(  t,x\right)  +h(t,x,g(t,x),\partial
_{x}g(t,x)\sigma(t,x,u^{\ast}))=0, &
\end{cases}
\label{hjb-71}%
\end{equation}
subject to $\tilde{v}(T,x)=x$ and $g(T,x)=\exp\left(  -\mathcal{\lambda
}x\right)  $. Note that, when $h(\cdot)=0$, $\mathbb{E}_{t}^{g}[\cdot]$
becomes a linear expectation. Then
\[
\tilde{Y}_{t}^{t,x,u}:=J(t,x;u(\cdot))
\]
is the solution of the following one dimensional quadratic BSDE:%
\begin{equation}
\tilde{Y}^{t,x,u}=X_{T}^{t,x,u}-\int_{t}^{T}\frac{\lambda(\tilde{Z}%
_{s}^{t,x,u})^{2}}{2}\mathrm{d}s+\int_{t}^{T}\tilde{Z}_{s}^{t,x,u}%
\mathrm{d}B_{s}.\label{quadratic}%
\end{equation}
Therefore, the HJB equation (\ref{hjb-71}) becomes the classical one
dimensional HJB equation, see, for instance, \citet{BLY19} and \citet{M20}.
{However, for the general case, $h(\cdot)\neq0$, \emph{problem (\ref{risksens}) falls outside the scope of the models developed by} \citet{Basak2010} or
\citet{BKM2017}. In this case, our extended HJB equation (\ref{hjb-71}) becomes indispensable for formulating the problem and identifying the corresponding equilibrium. Thus, our model extends
their formulations to a more
general framework. }

\section{{Conclusion}}\label{sec8}

In this paper, we establish an extended DPP for a general control problem in which the states satisfy a coupled FBSDE. Consequently, it leads to rigourous derivation of an extended HJB equation for the controlled system. We explore both smooth and viscosity solutions for this extended HJB equation.

Given that the value function (\ref{value-v}) relies on the solution of a controlled FBSDE, constructing an optimal feedback control $u^{\ast}(\cdot)$ becomes challenging. To the best of our knowledge, since \citet{P93}, previous studies have not delved into the DPP of the value function (\ref{value-v}) and its associated state-dependent PDE method. As the objective functional involves a BSDE state, the controlled system corresponds to an extended HJB equation, which manifests as a vector-valued PDE. Due to the scarcity of a comprehensive comparison theorem for vector-valued PDEs, the conventional viscosity method typically does not apply to such equations. Therefore, this paper introduces a novel definition of the viscosity solution, incorporating first-order smoothness in the auxiliary equation. Consequently, we establish the existence and uniqueness of the viscosity solution for the extended HJB system. Our results complement the existing literature on viscosity solutions to multidimensional PDEs.

Our method of extended DPP has the potential to provide equilibrium solutions for a wide array of time-inconsistent control problems. To illustrate its applicability, we examine problems related to the traditional mean-variance model, utility maximization for narrow framing investors and risk-sensitive control under a nonlinear expectation.

{In summary, in comparison to the existing maximum principle method, as exemplified in, for instance, \citet{P93}, this paper introduces a PDE-based approach to provide equilibrium solutions for a general control problem where the states are described by an FBSDE system, thereby filling a longstanding gap in the literature regarding PDE method for this class of control problems. Nevertheless, future research should explore further enhancements, such as addressing scenarios in which the FBSDEs are intricately coupled, with the forward diffusion process depending on the solution of the backward component. Notably, the initial value of the SDE depends on $(x+Y_{t})$. }

\end{document}